\documentclass{jcmlatex}
\usepackage{subfigure}
\usepackage{epstopdf}

\def\JCMvol{xx}
\def\JCMno{x}
\def\JCMyear{200x}
\def\JCMreceived{xxx}
\def\JCMrevised{xxx}
\def\JCMaccepted{xxx}

\usepackage{mathtools}

\allowdisplaybreaks

\usepackage{amssymb}
\usepackage[mathscr]{eucal}

\usepackage{mismath}

\usepackage{hyperref}

\usepackage{bm} 

\usepackage{graphicx}
\graphicspath{{figs/}}

\usepackage{tikz}
\usetikzlibrary{math}
\usetikzlibrary{calc}

\usepackage{subcaption}

\usepackage{multirow}

\usepackage{makecell}
\setcellgapes{2pt}

\usepackage{booktabs}

\usepackage{siunitx}
\usepackage{algorithm}
\usepackage{algpseudocode}

\newcommand{\RomanNumeralCaps}[1]
{\MakeUppercase{\romannumeral #1}}

\newcommand*\circled[1]{\tikz[baseline=(char.base)]{
    \node[shape=circle,draw,inner sep=1pt] (char) {\textbf{\small #1}};}}

\usepackage[skins]{tcolorbox}
\newtcolorbox{stepbox}[2][]{%
  enhanced,
  attach boxed title to top center={yshift=-3mm,yshifttext=-1mm},
  colframe=blue!75!black,
  colbacktitle=red!80!black,
  fonttitle=\bfseries,
  title=#2,#1
}

\newcommand{\SubplotTag}[1]{\textbf{\small \textsf{#1}}}

\DeclareMathOperator{\Div}{div}
\DeclareMathOperator{\Dist}{dist}
\DeclareMathOperator{\Span}{span}

\DeclareMathOperator{\Ess}{ess}
\DeclareMathOperator{\Supp}{supp}
\DeclareMathOperator{\Int}{int}
\DeclareMathOperator{\Diam}{diam}
\DeclareMathOperator{\Argmin}{argmin}

\DeclarePairedDelimiter{\RoundBrackets}{(}{)}
\DeclarePairedDelimiter{\CurlyBrackets}{\{}{\}}
\DeclarePairedDelimiter{\SquareBrackets}{[}{]}

\usepackage[capitalise]{cleveref}

\crefname{assumption}{assumption}{assumptions}
\Crefname{assumption}{Assumption}{Assumptions}

\begin{document}
\markboth{E. CHUNG, P. CIARLET, X. JIN AND C. YE}{Multiscale modeling for a class of high-contrast heterogeneous sign-changing problems}

\title{Multiscale modeling for a class of high-contrast heterogeneous sign-changing problems}

\author{Eric T. Chung
\thanks{Department of Mathematics,The Chinese University of Hong Kong, Hong Kong, China \\ Email: tschung@math.cuhk.edu.hk}
\and
 Patrick Ciarlet Jr.
\thanks{POEMS, CNRS, INRIA, ENSTA Paris, Institut Polytechnique de Paris, 91120 Palaiseau, France\\ Email: patrick.ciarlet@ensta.fr}
\and
Xingguang Jin \footnote{Corresponding author}
\thanks{Department of Mathematics,The Chinese University of Hong Kong, Hong Kong, China \\ Email: xgjin@math.cuhk.edu.hk }
\and
Changqing Ye
\thanks{Department of Mathematics,The Chinese University of Hong Kong, Hong Kong, China \\ Email: cqye@math.cuhk.edu.hk}}
\maketitle

\begin{abstract}
 The mathematical formulation of sign-changing problems involves a linear second-order partial differential equation in the divergence form, where the coefficient can assume positive and negative values in different subdomains.
These problems find their physical background in negative-index metamaterials, either as inclusions embedded into common materials as the matrix or vice versa.
In this paper, we propose a numerical method based on the constraint energy minimizing generalized multiscale finite element method (CEM-GMsFEM) specifically designed for sign-changing problems.
The construction of auxiliary spaces in the original CEM-GMsFEM is tailored to accommodate the sign-changing setting.
The numerical results demonstrate the effectiveness of the proposed method in handling sophisticated coefficient profiles and the robustness of coefficient contrast ratios.
Under several technical assumptions and by applying the \texttt{T}-coercivity theory, we establish the inf-sup stability and provide an a priori error estimate for the proposed method.
\end{abstract}



\section{Introduction}
Metamaterials greatly expand the design space of materials by offering unconventional properties that are not found in nature.
Typically, metamaterials are created by assembling common materials periodically.
The exotic characteristics of metamaterials manifest at the effective medium level, which can strikingly differ from the intrinsic properties of the constituent materials.
Some notable examples include auxetics, exhibiting a negative Poisson's ratio \cite{Almgern1985,Lakes2017}; pentamode materials, characterized by a vanishing shear modulus \cite{Milton1995,Kadic2012}; and negative-index materials, displaying negative electric and magnetic permeability \cite{Veselago1968,Shelby2001}.
From a mathematical perspective, the emergence of metamaterials prompts a reconsideration of the well-posedness theory of Partial Differential Equations (PDEs), as the coefficients in PDEs may fall outside the conventional range for coercivity.
In this context, we focus on the so-called sign-changing problem, which is rooted in the background that negative-index materials are embedded into a common medium, or vice versa.

The mathematical nature of the sign-changing problem is a linear second-order PDE in the divergence form, where the coefficient allows both positive and negative values in different subdomains, with a discontinuity across the interface between the subdomains.
The well-posedness of the model problem is generally questionable due to the absence of uniform strict positivity of the coefficient.
To address this, \texttt{T}-coercivity was introduced by Bonnet-BenDhia, Ciarlet Jr., and Zw\"{o}lf in their work \cite{BonnetBenDhia2010}, providing a systematic approach to study the well-posedness of the model problem.
The key idea behind \texttt{T}-coercivity is to devise a bijective map $\mathcal{T}$ that enables the inf-sup (Banach--N\v{e}cas--Babu\v{s}ka) condition to hold trivially.
The construction of $\mathcal{T}$ only relies on the geometric information of interfaces that separate the subdomains with different signs of the coefficient.
The main statement of the \texttt{T}-coercivity theory is that if the contrast ratio between the negative and positive coefficients is sufficiently large, then the model problem becomes well-posed.
In \cite{BonnetBenDhia2012}, the sharpest condition of the contrast ratio ensuring well-posedness on certain simple interface problems was derived.
Further extensions of the \texttt{T}-coercivity theory encompass Helmholtz-like problems \cite{CiarletJr.2012}, time-harmonic Maxwell equations \cite{BonnetBenDhia2014,BonnetBenDhia2014a}, eigenvalue problems \cite{Carvalho2017}, and mixed problems \cite{Barre2023}.

Although analytical solutions offer us a profound mathematical insight into sign-changing problems, numerical methods are indispensable for practical applications.
The \texttt{T}-coercivity theory justifies the well-posedness of the PDE models, which also guides the design of suitable approximations of the original problem.
An intriguing aspect lies in the compatibility between $\mathcal{T}$ at the continuous level and the discretization level.
As an immediate result, it is emphasized in \cite{Chesnel2013} that, to obtain an optimal convergence rate, meshes near flat interfaces should be symmetric.
Later, a new treatment at the corners of interfaces was proposed in \cite{BonnetBenDhia2018}, leading to meshing rules for an arbitrary polygonal interface.
Considering the low regularity of the solution due to the heterogeneous coefficient, a posteriori error analysis was conducted in \cite{Nicaise2011} and in \cite{CiarletJr.2018}, which provides a reliable error estimator for adaptive mesh refinement routines \cite{Verfuerth2013}.
Beyond classic continuous Galerkin methods, Chung and Ciarlet Jr.\ introduced a staggered discontinuous Galerkin method in \cite{Chung2013}, accompanied by stability and convergence analysis.

The nature of heterogeneous coefficients in sign-changing problems motivates the application of multiscale computational methods.
Pioneered by Hou and Wu in \cite{Hou1997}, the methodology of incorporating model information into the construction of finite element spaces, coined as MsFEMs, has garnered significant attention.
As a discretization scheme, an advantage of MsFEMs is that the meshes are not required to resolve the heterogeneity of the coefficient, although the implementation of MsFEMs indeed relies on a pair of nested meshes.
To relieve the rigidity from boundary conditions in constructing multiscale bases, the oversampling technique is introduced in \cite{Hou1997} and subsequently proved to improve convergence rates (ref.\ \cite{Efendiev2000,Efendiev2009}).
The accuracy of MsFEMs, to some extent, may deteriorate when the coefficient violates the scale-separation assumption, as evidenced by the convergence theories in \cite{Hou1999,Ye2020,Ming2024}.
To address this issue, the Generalized Multiscale Finite Element Method (GMsFEM) was proposed by Efendiev, Galvis, and Hou in \cite{Efendiev2013}.
GMsFEMs leverage spectral decomposition to perform dimension reduction for the online space, exhibiting superior performance when dealing with high-contrast and channel-like coefficient profiles \cite{Chung2014a}.
The first construction of multiscale bases capable of achieving the theoretically best approximation property for general $L^\infty$ coefficients was credited to M\aa{}lqvist and Peterseim in their celebrating work \cite{Maalqvist2014}.
This construction, known as Localized Orthogonal Decomposition (LOD), utilizes quasi-interpolation operators to decompose the solution into macroscopic and microscopic components \cite{Altmann2021,Maalqvist2021}.
The combination of GMsFEMs and LOD led to the development of a CEM-GMsFEM by Chung, Efendiev, and Leung in \cite{Chung2018}, where ``CEM'' is the acronym for ``Constraint Energy Minimizing''.
The novelty of CEM-GMsFEMs resides in replacing quasi-interpolation operators in LOD with element-wise eigenspace projections.
Moreover, CEM-GMsFEMs introduce a relaxed version of the energy minimization problems to construct multiscale bases, which eliminates the necessity of solving saddle-point linear systems.
Our intention here is not to present a comprehensive review of multiscale computational methods from the community, and hence, notable advancements such as heterogeneous multiscale methods \cite{Weinan2003, Abdulle2012}, generalized finite element methods \cite{Babuska2011,Babuska2020,Ma2022}, and variational multiscale methods \cite{Hughes1995,Hughes2007} are not covered.

This article serves as an application of the CEM-GMsFEM to sign-changing problems.
While the \texttt{T}-coercivity theory guarantees the well-posedness of the model problem, the heterogeneity of the coefficient leads to a generally low regularity of the solution, resulting in suboptimal convergence rates when using standard finite element methods.
The CEM-GMsFEM, being a multiscale computational method, is specifically designed to handle the low regularity of the solution, and the construction of multiscale bases in this paper is tailored to the sign-changing setting.
For instance, the auxiliary space forms a core module in the original CEM-GMsFEM and is created by solving generalized eigenvalue problems, where the coefficient enters bilinear forms on both sides.
However, this approach is not suitable for sign-changing problems, as the eigenvalues can be negative, rendering the generalized eigenvalue problems ill-defined.
To ensure the positivity of the eigenvalues and the approximation ability of the auxiliary space, we modify the generalized eigenvalue problems by replacing the coefficient with its absolute value counterpart.
It is worth noting that in building multiscale bases, we adhere to the relaxed version of the energy minimization problems in this paper, which offers implementation advantages.
We conduct numerical experiments to validate the performance of the proposed method, emphasizing that coarse meshes need not align with the sign-changing interfaces.
Moreover, we demonstrate that contrast robustness, which is an important feature of the original CEM-GMsFEM, is inherited in the proposed method.
Under several assumptions, we prove the existence of multiscale bases in the \texttt{T}-coercivity framework, the exponential decay property for the oversampling layers, and the inf-sup stability of the online space.
Moreover, we can provide an apriori estimate for the proposed method, indicating that errors can be bounded in terms of coarse mesh sizes and the number of oversampling layers, while being independent of the regularity of the solution.

Currently, to the best of our knowledge, efforts to apply multiscale computational methods to sign-changing problems are scarce.
An exception can be found in the work by Chaumont-Frelet and Verf\"{u}rth in \cite{ChaumontFrelet2021}, where they employed the LOD framework.
Note that contrast ratios play a significant role in the \texttt{T}-coercivity theory for assessing the well-posedness of the model problem.
In comparison, the proposed method presented in this paper is more general, as it can handle a wider range of coefficient profiles, including those with high contrast ratios.

This paper is organized as follows.
In \cref{sec:preliminaries}, we introduce the model problem and present the \texttt{T}-coercivity theory.
The construction of multiscale bases in the proposed method is detailed in \cref{sec:methods}.
To validate the performance of the proposed method, \cref{sec:numerical experiments} presents numerical experiments conducted on four different models.
All theoretical analysis for the proposed method is gathered in \cref{sec:analysis}.
Finally, in \cref{sec:conclusions}, we conclude the paper.

\section{Preliminaries}\label{sec:preliminaries}
For simplicity, we consider a 2D Lipschitz domain denoted as $\Omega$, which can be divided into two non-overlapping subdomains $\Omega^+$ and $\Omega^-$ with $\Gamma$ as the interfaces.
The extension of the following discussion to 3D is straightforward.
Let $\sigma$ belong to $L^\infty(\Omega)$ such that the essential infimum of $\sigma$ over $\Omega^+$ is greater than zero ($\Ess\inf_{x\in \Omega^+}\sigma(x)>0$), and the essential supremum of $\sigma$ over $\Omega^-$ is less than zero ($\Ess\sup_{x\in \Omega^-} \sigma(x)<0$).
In particular, $\sigma$ is discontinuous across the interface seperating $\Omega^+$ and $\Omega^-$.
We further introduce notations: $\sigma^+_{\mathup{max}}=\Ess\sup_{x\in \Omega^+} \sigma(x)$, $\sigma^+_{\mathup{min}}=\Ess\inf_{x\in \Omega^+} \sigma(x)$, $\sigma^-_{\mathup{max}}=\Ess\sup_{x\in \Omega^-} \abs{\sigma(x)}$ and $\sigma^-_{\mathup{min}}=\Ess\inf_{x\in \Omega^-} \abs{\sigma(x)}$.
Hence, we require that $\sigma_{\mathup{max}}^+ \geq \sigma_{\mathup{min}}^+ > 0$ and $\sigma_{\mathup{max}}^- \geq \sigma_{\mathup{min}}^- > 0$.
We redefine $V$ as the conventional Hilbert space $H^1_0(\Omega)$, and consider a bilinear form $a(\cdot, \cdot)$ on $V\times V$ given by:
\[
  a(v, w)=\int_\Omega \sigma \nabla v\cdot \nabla w \di x,\ \forall (v, w)\in V\times V.
\]
Then for an function $f$ belonging to $L^2(\Omega)$, the following variational form defines the model problem:
\begin{equation}\label{eq:model problem}
  \text{find } u\in V \text{ s.t. } \forall v \in V,\ a(u, v)=\int_\Omega f v \di x.
\end{equation}
Certainly, \cref{eq:model problem} corresponds to a PDE of $u$, i.e., $-\Div \sigma \nabla u = f$ with a boundary condition $u=0$ on $\partial \Omega$.
It is also convenient to introduce a notation for a norm as $\norm{v}_{\tilde{a},\omega}\coloneqq(\int_\omega \abs{\sigma} \abs{\nabla v}^2 \di x)^{1/2}$, where $\omega$ is a subdomain of $\Omega$.
Moreover, we abbreviate $\norm{\cdot}_{\tilde{a},\Omega}$ as $\norm{\cdot}_{\tilde{a}}$, which is exactly the \emph{energy norm} on $V$.

However, the well-posedness of \cref{eq:model problem} is generally questionable due to the loss of uniform strict positivity of $\sigma$.
\texttt{T}-coercivity is based on a fact: if there exists a \emph{bijective} map $\mathcal{T}\colon V\rightarrow V$ and a positive constant $\alpha$ such that for all $v\in V$, $\abs{a(v, \mathcal{T}v)} \geq \alpha \norm{v}_{\tilde{a}}^2$, then the Banach--N\v{e}cas--Babu\v{s}ka theory confirms the existence and uniqueness of the solution.
The novelty of \texttt{T}-coercivity lies in providing a systematic construction of such $\mathcal{T}$ using a ``flip'' operator, which is essentially dependent on the geometric information of the interfaces between $\Omega^+$ and $\Omega^-$.
Let $V^\pm \subset H^1(\Omega^\pm)$ be the space obtained by restricting functions in $V$ in $\Omega^\pm$ and $\Gamma$ be the interface.
Assuming the existence of a bounded map $\mathcal{R}\colon V^+\rightarrow V^-$ with $\mathcal{R}v|_{\Gamma}=v|_{\Gamma}$ in the sense of Sobolev traces, the construction of $\mathcal{T}$ is given by
\begin{equation}\label{eq:T operator}
  \mathcal{T}v=\begin{cases}
    v_1,                  & \text{ in } \Omega^+, \\
    -v_2+2\mathcal{R}v_1, & \text{ in } \Omega^-,
  \end{cases}
\end{equation}
where
\begin{equation}\label{eq:split of v}
  v=\begin{cases}
    v_1, & \text{ in } \Omega^+, \\
    v_2, & \text{ in } \Omega^-.
  \end{cases}
\end{equation}
Certainly, the definition of $\mathcal{T}$ in \cref{eq:T operator} is bijective and bounded. Furthermore, we can show that
\begin{equation}\label{eq:derivation of T-coercivity}
\begin{split}
a(v, \mathcal{T}v) & = \int_{\Omega^+}\abs{\sigma}\nabla v_1\cdot \nabla v_1 \di x - \int_{\Omega^-} \abs{\sigma} \nabla v_2 \cdot \nabla (-v_2 + 2\mathcal{R}v_1)\di x\\
& \underset{(\forall \eta > 0)}{\geq} \int_{\Omega^+}\abs{\sigma}\abs{\nabla v_1}^2 \di x-\frac{1}{\eta} \int_{\Omega^-} \abs{\sigma}\abs{\nabla \mathcal{R}v_1}^2\di x  + (1-\eta)\int_{\Omega^-}\abs{\sigma}\abs{\nabla v_2}^2\di x \quad \text{\footnotesize (by Young's inequality)}\\
& \geq \int_{\Omega^+}\abs{\sigma}\abs{\nabla v_1}^2 \di x - \frac{\sigma_{\mathup{max}}^- \norm{\mathcal{R}}_{1}^2}{\eta} \int_{\Omega^+} \abs{\nabla v_1}^2\di x +  (1-\eta)\int_{\Omega^-}\abs{\sigma}\abs{\nabla v_2}^2\di x \\
& \geq \RoundBrackets*{1-\frac{\sigma_{\mathup{max}}^-\norm{\mathcal{R}}_1^2}{ \sigma_{\mathup{min}}^+\eta}} \int_{\Omega^+}\abs{\sigma}\abs{\nabla v_1}^2 \di x + (1-\eta)\int_{\Omega^-}\abs{\sigma}\abs{\nabla v_2}^2\di x \\
& \geq \RoundBrackets*{1-\sqrt{\frac{\sigma_{\mathup{max}}^-}{\sigma_{\mathup{min}}^+}}\norm{\mathcal{R}}_1} \norm{v}_{\tilde{a}}^2, \quad \text{\footnotesize (by choosing $\eta=\sqrt{{\sigma_{\mathup{max}}^-}/{\sigma_{\mathup{min}}^+}}\norm{\mathcal{R}}_1$)}
\end{split}
\end{equation}
where $\norm{\mathcal{R}}_1$ should be understood as a positive constant such that for all $w\in V^+$,
\begin{equation}\label{eq:norm1 of R}
\RoundBrackets*{\int_{\Omega^-} \abs{\nabla \mathcal{R}w}^2 \di x}^{1/2} \leq \norm{\mathcal{R}}_1 \RoundBrackets*{\int_{\Omega^+} \abs{\nabla w}^2 \di x}^{1/2}.
\end{equation}
We can observe that \cref{eq:model problem} will always be well-posed if the ratio $\sigma_{\mathup{min}}^+/\sigma_{\mathup{max}}^-$ is large enough.
Note that the derivation in \cref{eq:derivation of T-coercivity} involves flipping ``positive'' to ``negative'' as indicated by the definition of  $\mathcal{T}\colon V^+\rightarrow V^-$.
It is possible to flip ``negative'' to ``positive'', and the well-posedness of \cref{eq:model problem} can also be achieved if the $\sigma_{\mathup{min}}^-/\sigma_{\mathup{max}}^+$ is large, as discussed in \cite{BonnetBenDhia2012}.

To gain insight into the operator $\mathcal{R}$, let us consider a square domain $\Omega=(0, 1)\times (0, 1)$ with $\Omega^+=(0, 1)\times (\gamma, 1)$ and $\Omega^-=(0, 1)\times(0, \gamma)$, where $0 < \gamma < 1$. For a $w \in V^+$, we can define $\mathcal{R}w$ as follows:
\[
  (\mathcal{R}w)(x_1, x_2) = w(x_1, 1-\tfrac{1-\gamma}{\gamma}x_2),\ \forall (x_1, x_2) \in \Omega^-.
\]
It is also straightforward to calculate that $\norm{\mathcal{R}}_1=\max\CurlyBrackets{\sqrt{(1-\gamma)/{\gamma}}, 1}$.
Note that the construction of $\mathcal{R}$ is not unique. As an alternative approach, we can introduce a smooth cut-off function $\xi$ satisfying $\xi=1$ on $\Gamma$ and $\Supp \xi w \subset (0, 1) \times (\gamma, a)$ with $ \gamma < a <1$.
Then, we can utilize the same technique to design $\mathcal{R}'$ such that $\Supp \mathcal{R}'(\xi w) \subset (0, 1)\times(b, \gamma)$ with $0 \leq b < \gamma$.
Naturally, as an operator, $\mathcal{R}'\xi\colon V^+\rightarrow V^-$ is also valid.

\section{Methods} \label{sec:methods}
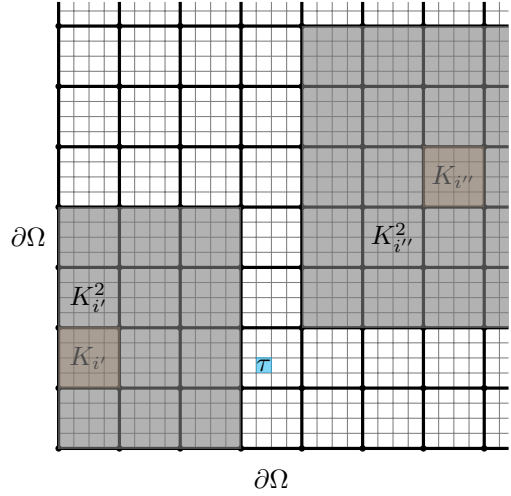
\begin{figure}[!ht]
  \centering
  \begin{tikzpicture}[scale=0.8]
    \draw[step=0.25, gray, thin] (0.0, 0.0) grid (7.4, 7.4);
    \draw[step=1.0, black, very thick] (0.0, -0.0) grid (7.4, 7.4);
    \foreach \x in {0,...,7}
    \foreach \y in {0,...,7}{
        \fill (1.0 * \x, 1.0 * \y) circle (1.5pt);
      }
    \fill[brown, opacity=0.4] (0.0, 1.0) rectangle (1.0, 2.0);
    \node at (0.5, 1.5) {$K_{i'}$};
    \fill[opacity=0.6, gray] (0.0, 0.0) rectangle (3.0, 4.0);
    \node at (0.5, 2.5) {$K_{i'}^2$};

    \fill[brown, opacity=0.4] (6.0, 4.0) rectangle (7.0, 5.0);
    \node at (6.5, 4.5) {$K_{i''}$};
    \fill[opacity=0.6, gray] (4.0, 2.0) rectangle (7.4, 7.0);
    \node at (5.5, 3.5) {$K_{i''}^2$};

    \fill [cyan, opacity=0.5] (3.25, 1.25) rectangle (3.5, 1.5);
    \node at (3.375, 1.375) {$e$};

    \node at (-0.5, 3.5) {$\partial \Omega$};
    \node at (3.5, -0.5) {$\partial \Omega$};
  \end{tikzpicture}
  \caption{
    Illustration of the nested meshes $\mathscr{K}_h$ and $\mathscr{K}_H$.
    A fine element $e$, two coarse elements $K_{i'}$ and $K_{i''}$, accompanied by their corresponding oversampling regions $K_{i'}^2$ and $K_{i''}^2$, are colored differently.
  }
  \label{fig:grid}
\end{figure}

The proposed method relies on a pair of nested meshes $\mathscr{K}_h$ and $\mathscr{K}_H$. The fine mesh $\mathscr{K}_h$ should be capable of resolving the heterogeneity of $\sigma$.
Meanwhile, the construction of multiscale bases, known as the offline phase, is performed on $\mathscr{K}_h$.
Solving the final linear system is commonly referred to as the online phase, whose computational complexity is associated with the coarse mesh $\mathscr{K}_H$.
We reserve the index $i$ for enumerating all coarse elements on $\mathscr{K}_H$, where $1 \leq i \leq N_\mathup{elem}$.
The oversampling region $K_i^m$ of a coarse element $K_i$ holds a vital significance in describing the method.
Mathematically, $K_i^m$ is defined by a recursive relation:
\[
  \forall m \geq 1,\quad
  K_i^m \coloneqq \Int\RoundBrackets*{\cup\CurlyBrackets{\overline{K}\mid K \in \mathscr{K}_H \text{ with } \overline{K} \cap \overline{K_i^{m-1}} \neq \varnothing}}
\]
and $K_i^0 \coloneqq K_i$.
\cref{fig:grid} serves as an illustration of the nested meshes, and note that two oversampling regions $K_{i'}^2$ and $K_{i''}^2$ are colored in gray.

In the original CEM-GMsFEM, a generalized eigenvalue problem,
\begin{equation}\label{eq:old aux}
  \text{find } \lambda \in \mathbb{R} \text{ and } v \in H^1(K_i)\setminus\CurlyBrackets{0} \text{ s.t. }
  \int_{K_i} \sigma \nabla v \cdot \nabla w \di x = \lambda \int_{K_i} \mu_\mathup{msh} \Diam(K_i)^{-2} \sigma v w \di x,
\end{equation}
is solved on each coarse element $K_i$, where $\mu_\mathup{msh}$ is a generic positive constant that depends on the mesh quality.
Then the local auxiliary space $V_i^\mathup{aux}$ is formed by collecting leading eigenvectors.
However, recalling that $\sigma$ is not uniformly positive, the left-hand bilinear form in \cref{eq:old aux} is not positive semidefinite, and similarly, the right-hand bilinear form in \cref{eq:old aux} is not positive definite.
Then, determining leading eigenvectors via \cref{eq:old aux} is problematic since the eigenvalues could be negative.
To address this issue, we instead construct the following generalized eigenvalue problem on each $K_i$, which forms the first step of the proposed method:
\begin{stepbox}{First step}
  \begin{equation}\label{eq:new aux}
    \begin{aligned}
       & \text{find } \lambda \in \mathbb{R} \text{ and } v \in H^1(K_i)\setminus\CurlyBrackets{0} \text{ s.t. } \forall w \in H^1(K_i),     \\
       & \int_{K_i} \abs{\sigma} \nabla v \cdot \nabla w \di x = \lambda \int_{K_i} \mu_\mathup{msh} \Diam(K_i)^{-2} \abs{\sigma} v w \di x.
    \end{aligned}
  \end{equation}
\end{stepbox}
One may argue that if $\sigma$ only takes values in $\{-1, 1\}$, then \cref{eq:new aux} would yield trivial eigenspaces corresponding to the Laplace operator, thus failing to capture any heterogeneity information of $\sigma$.
However, according to \texttt{T}-coercivity, in such cases, the well-posedness of the model problem is violated (ref.\ \cite{BonnetBenDhia2012} section 6).
This inherent issue suggests that the original model problem should require a more sophisticated investigation in this case, which falls outside the scope of the proposed method.
Upon solving the $l^*$ leading eigenvectors, denoted as $\psi_{i,j}$ with $1 \leq j \leq l^*$, we can construct the local auxiliary space $V_i^\mathup{aux} \subset L^2(K_i)$ as $\Span\{\psi_{i,j}\}$ for $1 \leq i \leq N_{\mathup{elem}}$.
The global auxiliary space $V^\mathup{aux} \subset L^2(\Omega)$ is defined as $V^\mathup{aux}=\oplus_{i=1}^{N_\mathup{elem}} \widetilde{V}_i^\mathup{aux}$, where $\widetilde{V}_i^\mathup{aux}\subset L^2(\Omega)$ is the space by performing zero-extension of functions in $V_i^\mathup{aux}$.

We introduce several notations for future reference. Let $\mu$ be a function in $L^\infty(\Omega)$ satisfying
\[
  \mu|_{K_i}=\mu_\mathup{msh}\Diam(K_i)^{-2}\sigma|_{K_i}
\]
for all $K_i \in \mathscr{K}_H$.
Next, for a subdomain $\omega$ in $\Omega$, we define two bilinear forms:
\[
  a(v, w)_{\omega} \coloneqq \int_{\omega} \sigma \nabla v \cdot \nabla w \di x \text{ and } s(v, w)_{\omega} \coloneqq \int_{\omega} \mu v w \di x.
\]
Similarly to the definition of $\norm{\cdot}_{\tilde{a},\omega}$, we define $\norm{v}_{\tilde{s},\omega}$ as $(\int_{\omega} \abs{\mu} \abs{v}^2 \di x)^{1/2}$.
Again, we drop the subscript $\omega$ if $\omega=\Omega$.
Additionally, we define the orthogonal projection operator $\mathcal{P}_H\colon L^2(\Omega) \rightarrow L^2(\Omega)$ under the norm $\norm{\cdot}{\tilde{s}}$, with $V^\mathup{aux}=\im \mathcal{P}_H$.
We also use the shorthand notation $V_i^m$ as $H^1_0(K_i^m)$, and we always implicitly identify $V_i^m$ as a subspace of $V$.

Functions in $V^\mathup{aux}$ may not be continuous in $\Omega$, and thus $V^\mathup{aux}$ cannot be used as a conforming finite element space.
In the second step, we will construct a multiscale basis $\phi_{i,j}$ in $V_i^m$, corresponding to $\psi_{i,j}$ in $V_i^\mathup{aux}$ that is obtained in the first step.
This construction is based on the following variational problem, where the righthand bilinear form is defined on the original coarse element $K_i$:
\begin{stepbox}{Second step}
  \begin{equation}\label{eq:basis}
    \begin{aligned}
       & \text{find } \phi_{i,j}\in V_i^m \text{ s.t. } \forall w \in V_i^m,                                                    \\
       & a(\phi_{i,j}, w)_{K_i^m} + s(\mathcal{P}_H\phi_{i,j}, \mathcal{P}_H w)_{K_i^m} = s(\psi_{i,j}, \mathcal{P}_H w)_{K_i}.
    \end{aligned}
  \end{equation}
\end{stepbox}
In the original CEM-GMsFEM \cite{Chung2018}, the variational form in \cref{eq:basis} is derived from a ``relaxed'' constrained energy minimization problem:
\[
  \phi_{i,j}=\Argmin\CurlyBrackets*{a(w, w)_{K_i^m}\mid \ w\in V_i^m,\ \mathcal{P}_H w = \psi_{i,j}}.
\]
We note that, in the present situation, the bilinear forms $a(\cdot, \cdot)_{K_i^m}$ and $s(\mathcal{P}_H\cdot, \mathcal{P}_H\cdot)_{K_i^m}$ may lack coercivity, leading to an ill-defined minimization.
The multiscale space is formed as
\[
  V_H^m = \Span\CurlyBrackets*{\phi_{i,j}\mid 1\leq i \leq N_\mathup{elem},\, 0 \leq j < l^*},
\]
and the solution of the model problem is approximated by solving the following variational problem:
\begin{equation}\label{eq:online}
  \text{find } u_H \in V_H^m \text{ s.t. } \forall w_H \in V_H^m,\ a(u_H, w_H)=\int_\Omega f w_H \di x.
\end{equation}
Note that \cref{eq:new aux,eq:basis} are practically solved on the fine mesh $\mathscr{K}_h$.

\begin{figure}[!ht]
  \includegraphics[width=\textwidth]{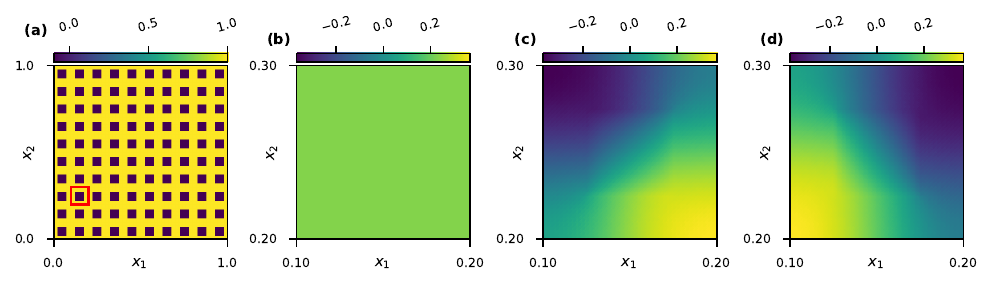}
  \caption{
    \SubplotTag{(a)} The coefficient profile and the marked coarse element.
    \SubplotTag{(b)}--\SubplotTag{(d)} The plot of the first/second/third eigenfunction corresponding to the marked coarse element.
  }\label{fig:eigen-0d1+1d0}
\end{figure}

\begin{figure}[!ht]
  \includegraphics[width=\textwidth]{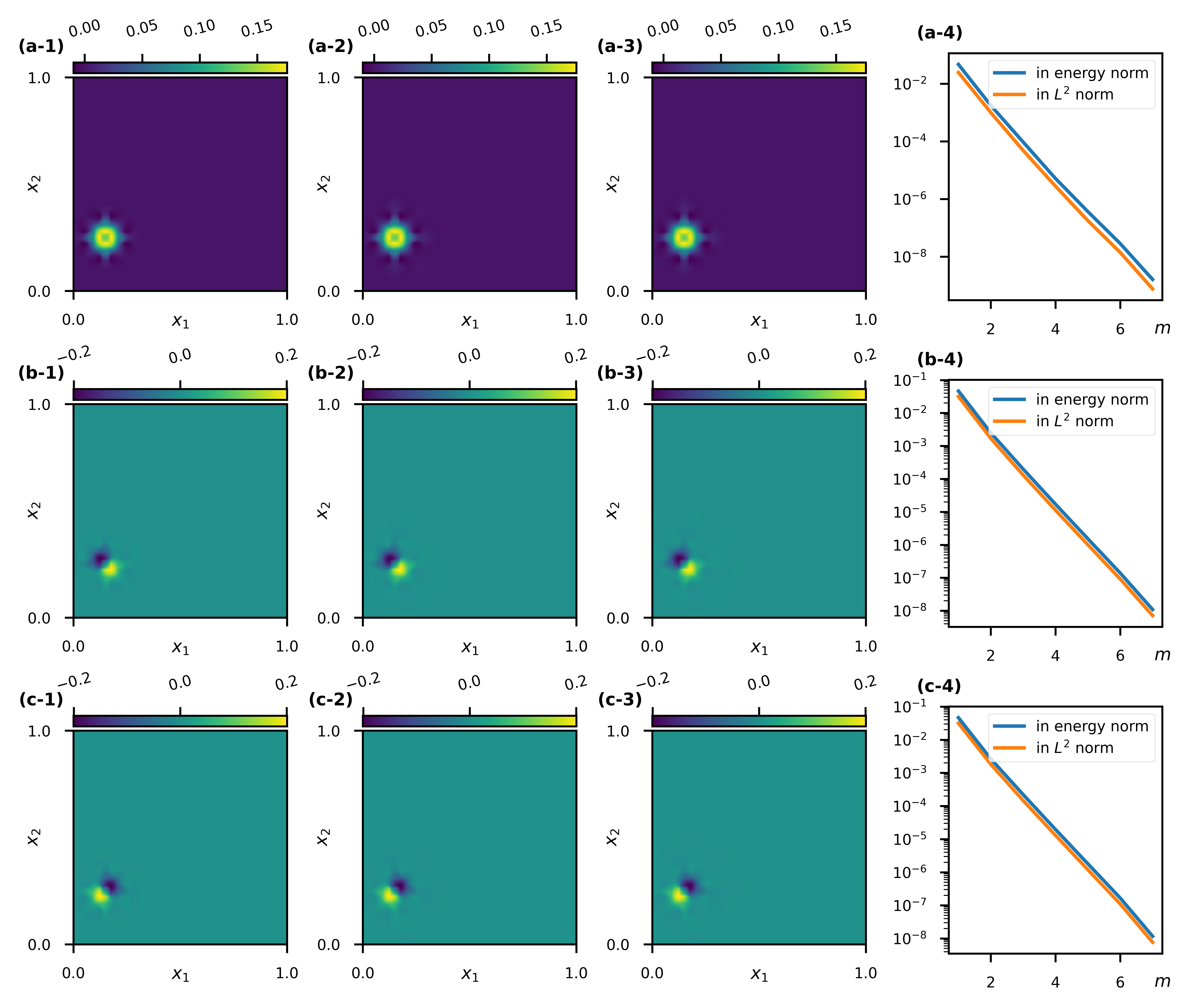}
  \caption{
    The subplots are marked as \SubplotTag{(x-y)}, where \SubplotTag{x} can take \SubplotTag{a}, \SubplotTag{b}, or \SubplotTag{c}, corresponding to the results for the first, second, or third eigenfunction, respectively.
    If \SubplotTag{y} is \SubplotTag{1}, \SubplotTag{2}, or \SubplotTag{3}, the subplot displays the multiscale basis with $m$ oversampling layers, $m$ equal to \SubplotTag{y}.
    Alternatively, if \SubplotTag{y} is \SubplotTag{4}, the subplot shows the relative differences (y-axis) in the energy and $L^2$ norm of the multiscale bases between $m=8$ and $m=1,\dots,7$ (x-axis).
  }\label{fig:ms-0d1+1d0}
\end{figure}

The success of the proposed method hinges on that the multiscale basis decays rapidly w.r.t.\ $m$, the number of oversampling layers.
To demonstrate this, we consider a square domain with a $10\times 10$ periodic structure, as illustrated in the subplot \SubplotTag{(a)} of \cref{fig:eigen-0d1+1d0}.
Each periodic cell contains a square inclusion that is centered with a negative coefficient imposed, resulting in the union of all inclusions forming the subdomain $\Omega^-$.
The length ratio of the inclusion to the periodic cell is set to $1/2$.
The material properties are set as $\sigma=1.0$ in $\Omega^+$ and $\sigma=-0.1$ in $\Omega^-$ such that \texttt{T}-coercivity is satisfied.
The coarse mesh $\mathscr{K}_H$ aligns with the periodic structure.
We select a coarse element marked with red borders in the subplot \SubplotTag{(a)} of \cref{fig:eigen-0d1+1d0} and plot the first three eigenfunctions $\psi_{i,1}$, $\psi_{i,2}$ and $\psi_{i,3}$ calculated via \cref{eq:new aux} in the subplots \SubplotTag{(b)} to \SubplotTag{(d)}.
It is worth noting that the first eigenvalue is always $0$, and the corresponding eigenfunction is always a constant function, as easily derived in \cref{eq:new aux} and also validated in the subplot \SubplotTag{(b)}.
The decay of the multiscale basis $\phi_{i,1}$, $\phi_{i,2}$ and $\phi_{i,3}$ are solved by \cref{eq:basis} with different oversampling layers is demonstrated in \cref{fig:ms-0d1+1d0}, while $a(\cdot,\cdot)_{K_i^m}$ and $s(\mathcal{P}_H \cdot, \mathcal{P}_H\cdot)_{K_i^m}$ are not coercive.
In this figure, the first, second, and third rows correspond to the results of the first, second, and third eigenfunctions, respectively, while the first, second, and third columns display the plots of the multiscale basis with $m=1$, $2$, and $3$, respectively.
The position of the selected coarse element determines the maximum value of $m$, which in this case is $8$.
Consequently, we calculate the relative differences in both the energy and $L^2$ norm of the multiscale bases between $m=8$ and $m=1,\dots,7$, and present the results in the fourth column of \cref{fig:ms-0d1+1d0}.
From the plots of multiscale bases, we observe that multiscale bases vanish away from the selected coarse element.
Hence, one may expect to accurately compute multiscale bases with a small number of oversampling layers.
Moreover, the relative differences in the logarithmic scale are almost linear w.r.t.\ $m$, which suggests the exponential decay may still hold in sign-changing problems.

\section{Numerical experiments}\label{sec:numerical experiments}
We conduct numerical experiments on a square domain $\Omega = (0,1) \times (0,1)$. The fine mesh $\mathscr{K}_h$ is generated by dividing $\Omega$ into $400 \times 400$ squares.
Consequently, the coefficient profile $\sigma$ is represented by a $400 \times 400$ matrix/image, with each element corresponding to a constant value on a fine element.
To investigate the convergence behavior of the proposed method with different coarse mesh sizes $H$, we consider four different coarse meshes $\mathscr{K}_H$: $10 \times 10$, $20 \times 20$, $40 \times 40$, and $80 \times 80$.
The reference solution, auxiliary spaces, and multiscale bases are all calculated on the fine mesh $\mathscr{K}_h$ using the $Q_1$ finite element method.
We evaluate the convergence of the proposed method using two relative error indices: the relative energy error and the relative $L^2$ error which are defined as follows:
\[
  \frac{\norm{e_h}_{\tilde{a}}}{\norm{u_h}_{\tilde{a}}} \text{  and  } \frac{\norm{e_h}_{L^2(\Omega)}}{\norm{u_h}_{L^2(\Omega)}},
\]
where $u_h$ is the reference solution calculated by the $Q_1$ FEM on $\mathscr{K}_h$ or the nodal interpolation of the exact solution (if available), and $e_h$ is the error between the reference solution and the numerical solution.
For simplicity, we take $\mu$ as
\[
  \mu|_{K_i}=\mu_\mathup{msh}\Diam(K_i)^{-2}\sigma|_{K_i}=24H^{-2}\sigma|_{K_i}
\]
for all numerical experiments, as suggested in \cite{Ye2023a}.
In the following discussions, we mark statements of direct observations from numerical experiments with a circled number, e.g., \circled{1}.
We implement the method using the Python libraries NumPy and SciPy\footnote{Instead of using the default sparse linear system solver in SciPy, we opted to utilize the \texttt{pardiso} solver to enhance efficiency.}, and all the codes are hosted on Github\footnote{\url{https://github.com/Laphet/sign-changing}}.

\subsection{Flat interface model}
We first consider a flat interface model described in the ending part of \cref{sec:preliminaries}, i.e., we define $\Omega^+=(0, 1)\times (\gamma, 1)$ and $\Omega^-=(0, 1)\times(0, \gamma)$ with $0 < \gamma < 1$.
In $\Omega^+$, we assign a fixed value of $\sigma^+_*$ to $\sigma^+$, while in $\Omega^-$, we take $-\sigma^-_*$ for $\sigma^-$, where $\sigma^+_*$ and $\sigma^-_*$ are both positive.
We can devise an exact solution $u$ as
\[
  u(x_1, x_2) = \begin{cases}
    -\sigma^-_*x_1(x_1-1)x_2(x_2-1)(x_2-\gamma), & \text{ in } \Omega^+, \\
    \sigma^+_*x_1(x_1-1)x_2(x_2-1)(x_2-\gamma),  & \text{ in } \Omega^-,
  \end{cases}
\]
which corresponds to a smooth source term $f$ given by
\[
  f(x_1, x_2) = \sigma^-_*\sigma^+_*\Big(2x_2(x_2 - 1)(x_2 - \gamma)+ x_1(x_1 - 1)(6x_2 - 2(\gamma + 1))\Big).
\]
According to the \texttt{T}-coercivity theory \cite{Chesnel2013}, the problem is well-posed if
\[
  \frac{\sigma^-_*}{\sigma^+_*} \notin \SquareBrackets*{\frac{\gamma}{1-\gamma}, 1},
\]
provided that $\gamma \leq 1/2$.

\begin{figure}[!ht]
  \centering
  \includegraphics[width=\textwidth]{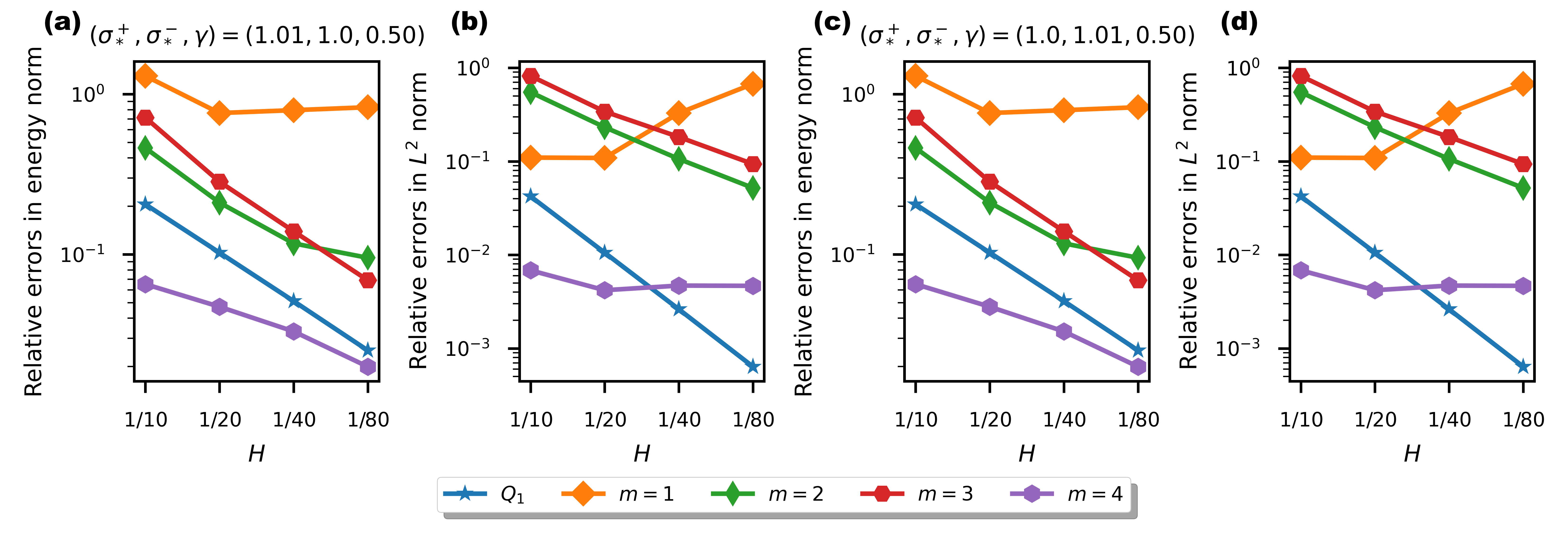}
  \caption{
    Numerical results for the flat interface model with $\gamma=1/2$, where the relative errors of the proposed method with different numbers of oversampling layers $m$ and the $Q_1$ FEM are calculated w.r.t.\ the coarse mesh size $H$.
    Subplots \SubplotTag{(a)} and \SubplotTag{(b)} correspond to $(\sigma^+_*, \sigma^-_*)=(1.01, 1)$, which the relative errors are measured in the energy and $L^2$ norm, respectively. Similarly, subplots \SubplotTag{(c)} and \SubplotTag{(d)} correspond to the setting $(\sigma^+_*, \sigma^-_*)=(1, 1.01)$, following the same manner.
  }\label{fig:flat-interface-error-l0-ps}
\end{figure}

\paragraph{Case \RomanNumeralCaps{1}} We set $\gamma=0.5$, we can check that now the interface $\Gamma = (0,1)\times \{\gamma\}$ is \emph{fully resolved by every coarse mesh}.
Therefore, we can expect that applying the $Q_1$ FEM on $\mathscr{K}_H$ directly can yield satisfactory accuracy. we conducted two groups of experiments with $(\sigma^+_*, \sigma^-_*)=(1.01, 1)$ and $(\sigma^+_*, \sigma^-_*)=(1, 1.01)$, both satisfying the \texttt{T}-coercivity condition.
For the setting of the proposed multiscale method, we fix $l^*=3$, indicating that we calculate the first three eigenfunctions in \cref{eq:new aux} and construct three multiscale bases for each coarse element, while we vary the oversampling layers $m$ from $1$ to $4$.
We refer to \cref{fig:flat-interface-error-l0-ps} for the numerical results.
\circled{1} We can observe from subplots \SubplotTag{(a)} to \SubplotTag{(d)} that the convergence of the $Q_1$ FEM manifests a linear pattern w.r.t.\ $H$ in the logarithmic scale, consistent with the theoretical expectation.
\circled{2} We can also see that the number of oversampling layers $m$ has a significant impact on the accuracy of the proposed method.
\circled{3} However, for the same $m$, the error decaying w.r.t.\ $H$ does not always hold, as depicted in subplots \SubplotTag{(b)} and \SubplotTag{(d)}.
\circled{4} Although, for $m=4$, the proposed method exhibits higher accuracy than the $Q_1$ FEM, the computational cost is significantly higher due to the sophisticated process of constructing multiscale bases.
Therefore, the proposed method seems more suitable for scenarios involving intricate coefficient profiles.
Interestingly, we notice that subplots \SubplotTag{(a)} and \SubplotTag{(c)}, as well as \SubplotTag{(b)} and \SubplotTag{(d)}, are almost identical, implying that the contrast ratio $\sigma^-_*/\sigma^+_*$ crossing the critical value $1$ does not generate a significant influence on numerical methods.

\begin{figure}[!ht]
  \centering
  \includegraphics[width=\textwidth]{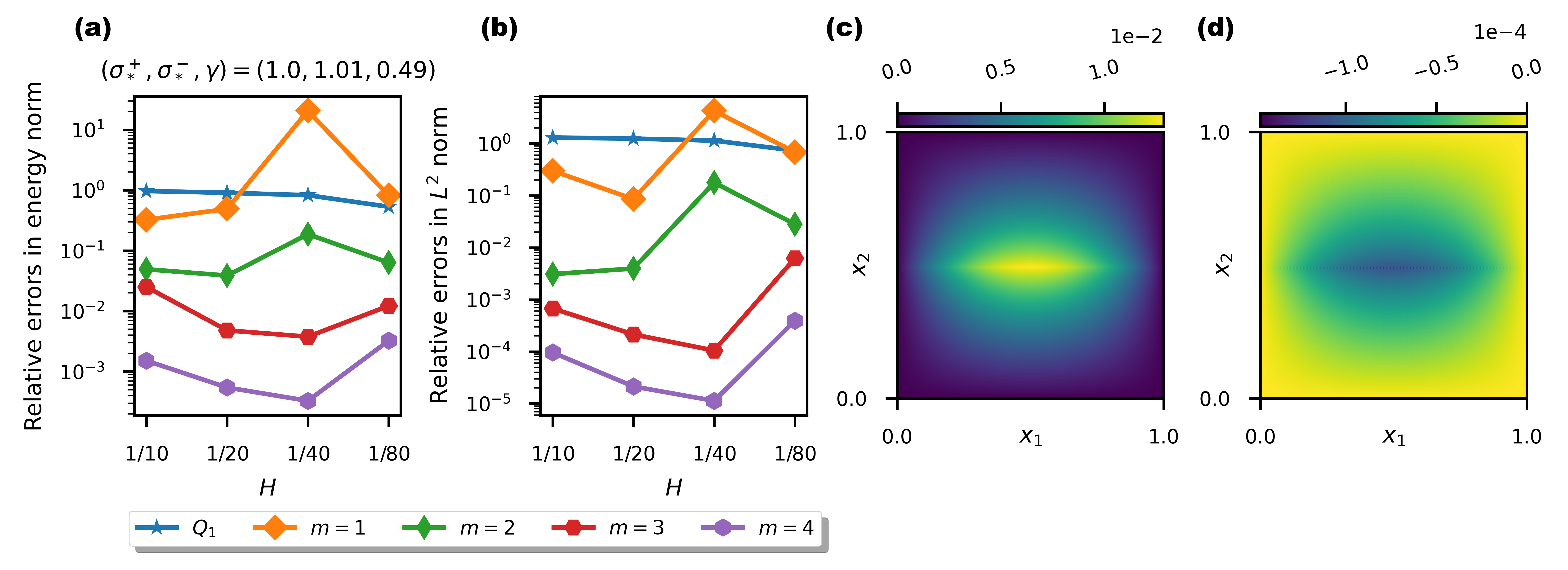}
  \caption{
    Numerical results for the flat interface model with $(\sigma^+_*, \sigma^-_*, \gamma)=(1, 1.01, 0.49)$.
    Subplots \SubplotTag{(a)} and \SubplotTag{(b)} show the relative errors of the proposed method with different numbers of oversampling layers $m$ and the $Q_1$ FEM w.r.t.\ the coarse mesh size $H$, but measured in different norms.
    Subplots \SubplotTag{(c)} and \SubplotTag{(d)} display the actual pointwise differences between the reference solution and numerical solutions obtained by the $Q_1$ FEM with $H=1/80$ and the proposed method with $(H, m)=(1/80, 3)$.
  }\label{fig:flat-interface-error-l1-ps}
\end{figure}

\paragraph{Case \RomanNumeralCaps{2}} Next, we consider the case that $\gamma = 0.49$, resulting in \emph{none of the coarse meshes being capable of resolving the interface}.
The results by setting $(\sigma^+_*, \sigma^-_*)=(1, 1.01)$ are reported in \cref{fig:flat-interface-error-l1-ps}, where subplots \SubplotTag{(a)} and \SubplotTag{(b)} correspond to the relative errors, measured in the energy and $L^2$ norm respectively.
Subplots \SubplotTag{(c)} and \SubplotTag{(d)} display the actual differences between the reference solution and numerical solutions obtained by the $Q_1$ FEM with $H=1/80$ and the proposed method with $(H, m)=(1/80, 3)$, respectively.
\circled{1} We observe that a slight change in $\gamma$ from $0.5$ to $0.49$ disrupts the convergence of the $Q_1$ FEM.
\circled{2} In contrast, the proposed method can achieve satisfactory accuracy at a level of $1\%$ for $m=3$ and of $0.1\%$ for $m=4$ in the energy norm.
\circled{3} Again, for a fixed $m$, the errors from the proposed method do not always decay w.r.t.\ $H$.
As shown in subplots \SubplotTag{(c)} and \SubplotTag{(d)}, both methods exhibit a concentration of errors near the interface, while the proposed method outperforms the $Q_1$ FEM by an order of magnitude of two in terms of pointwise errors.

\subsection{Periodic square inclusion model}
\begin{figure}[!ht]
  \centering
  \includegraphics[width=0.75\textwidth]{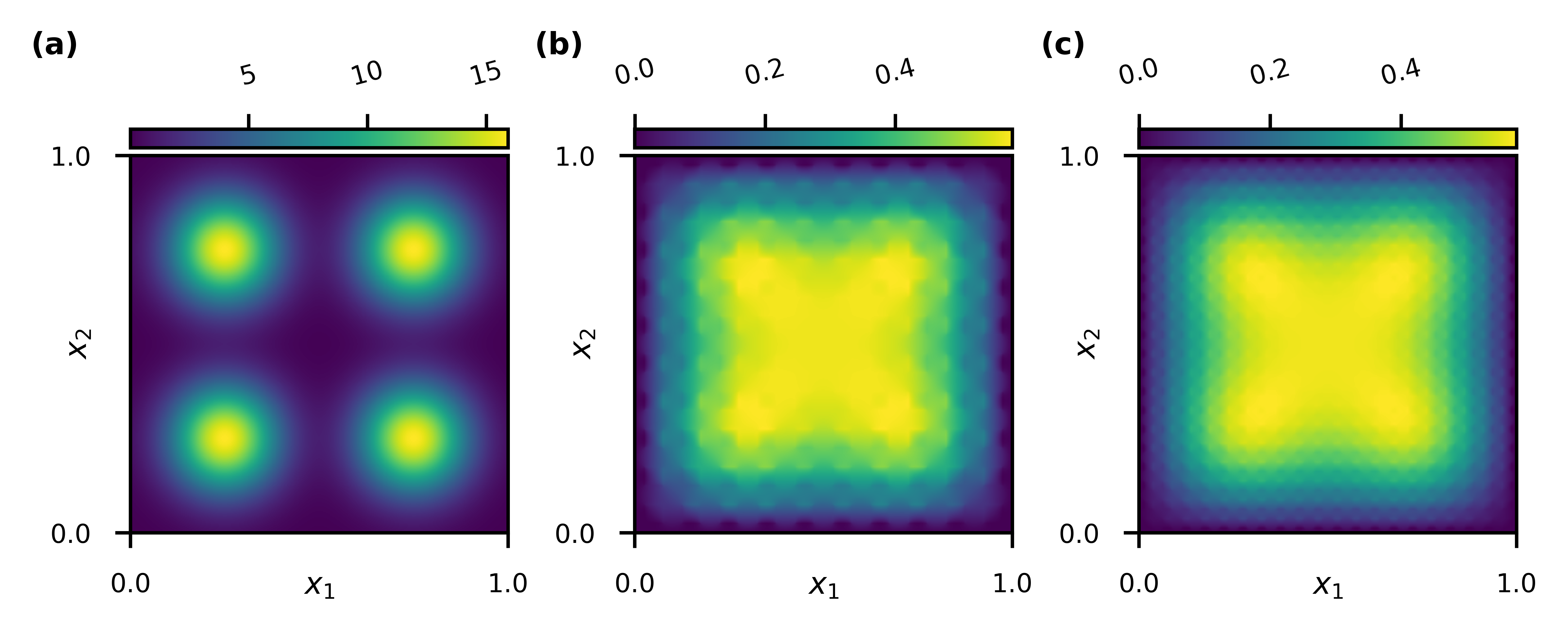}
  \caption{
    \SubplotTag{(a)} The smooth source term $f$ that is constructed as the superposition of four 2D Gaussian functions centered at $(0.25, 0.25)$, $(0.75, 0.25)$, $(0.25, 0.75)$, and $(0.75, 0.75)$ with a variance of $0.01$.
    \SubplotTag{(b)-(c)} The reference solutions correspond to the $10\times 10$ and $20\times 20$ periodic configurations, respectively.
  }\label{fig:square-inclusion-source-u}
\end{figure}

In this subsection, we revisit the periodic square inclusion model described in the ending part of \cref{sec:methods}.
The coefficient profile $\sigma$ is determined as \cref{fig:eigen-0d1+1d0}-\SubplotTag{(a)}.
Besides the $10\times 10$ periodic configuration shown in \cref{fig:eigen-0d1+1d0}-\SubplotTag{(a)}, we also consider the $20\times 20$ periodic configurations.
The source term $f$ is constructed as the superposition of four 2D Gaussian functions centered at $(0.25, 0.25)$, $(0.75, 0.25)$, $(0.25, 0.75)$, and $(0.75, 0.75)$ with a variance of $0.01$, as shown in \cref{fig:square-inclusion-source-u}-\SubplotTag{(a)}.
The reference solutions corresponding to the $10\times 10$ and $20\times 20$ periodic configurations are plotted in \cref{fig:square-inclusion-source-u}-\SubplotTag{(b)} and \SubplotTag{(c)}, respectively. We can observe that multiscale features emerge in the reference solutions, which may inspire future investigations into extending classical multiscale asymptotic analysis (ref.\ \cite{Bensoussan2011,Cioranescu1999}) to sign-changing problems.
Notably, the homogenization theory for the model was recently completed by Bunoiu and Ramdani in \cite{Bunoiu2016,Bunoiu2021,Bunoiu2022,Bunoiu2023}.
Meanwhile, numerical methods such as LOD \cite{ChaumontFrelet2021}, along with the proposed method, seek a low-dimensional representation of the solutions and can be regarded as general numerical homogenization techniques beyond the periodic setting.

We conduct numerical experiments using the proposed method on the $10\times 10$ and $20\times 20$ periodic configurations.
The relative errors in the energy norm and the $L^2$ norm are tabulated in \cref{tab:square-inclusion-error-10,tab:square-inclusion-error-20}, respectively.
We again take $l^*=3$ for all tests, while varying the number of oversampling layers $m$ from $1$ to $4$.
\circled{1} From \cref{tab:square-inclusion-error-10,tab:square-inclusion-error-20}, when $H=1/10$, enlarging $m$ from $2$ to $4$ does not lead to a significant improvement in accuracy, and relative errors in the energy norm reach a saturation level of $1\%$.
\circled{2} In comparison, for $H=1/80$, the numerical solutions achieve a relative error of $0.1\%$ in the energy norm for $m=3$ and $4$.
We can conclude that the accuracy of the proposed method is mutually controlled by the number of oversampling layers $m$ and the coarse mesh size $H$, and this convergence pattern is consistent with the flat interface model.
\circled{3} Another noteworthy observation is that the scale of models, i.e., $1/10$ and $1/20$, does not affect the accuracy of the proposed method, as the relative errors in \cref{tab:square-inclusion-error-10,tab:square-inclusion-error-20} are almost at the same level with different $H$ and $m$.
Hence, we can infer that the proposed method is robust to the scale of models and is not subject to the resonance error phenomenon encountered in the classic MsFEMs \cite{Hou1997,Hou1999,Efendiev2009}.

\begin{table}[!ht]
  \caption{
  For the $10\times 10$ periodic square inclusion model, the relative errors in the energy norm (in the columns labeled with $\norm{\cdot}_{\tilde{a}}$) and in the $L^2$ norm (in the columns labeled with $\norm{\cdot}_{L^2(\Omega)}$).
  }\label{tab:square-inclusion-error-10}
  \centering
  \makegapedcells
  \footnotesize{
    \begin{tabular}{c c c c c c c c c}
      \toprule
      \multirow{3}{*}{$H$} & \multicolumn{2}{c}{$m=1$}  & \multicolumn{2}{c}{$m=2$}    & \multicolumn{2}{c}{$m=3$}  & \multicolumn{2}{c}{$m=4$}                                                                                                                            \\
      \cmidrule{2-9}
                           & $\norm{\cdot}_{\tilde{a}}$ & $\norm{\cdot}_{L^2(\Omega)}$ & $\norm{\cdot}_{\tilde{a}}$ & $\norm{\cdot}_{L^2(\Omega)}$ & $\norm{\cdot}_{\tilde{a}}$ & $\norm{\cdot}_{L^2(\Omega)}$ & $\norm{\cdot}_{\tilde{a}}$ & $\norm{\cdot}_{L^2(\Omega)}$ \\
      \midrule
      $\frac{1}{10}$       & \num{2.433e-01}            & \num{8.931e-02}              & \num{5.162e-02}            & \num{5.741e-03}              & \num{5.225e-02}            & \num{5.785e-03}              & \num{5.232e-02}            & \num{5.786e-03}              \\
      $\frac{1}{20}$       & \num{3.785e-01}            & \num{1.923e-01}              & \num{4.960e-02}            & \num{4.981e-03}              & \num{5.583e-02}            & \num{5.922e-03}              & \num{5.662e-02}            & \num{6.024e-03}              \\
      $\frac{1}{40}$       & \num{6.978e-01}            & \num{5.833e-01}              & \num{4.011e-02}            & \num{2.292e-03}              & \num{1.753e-03}            & \num{3.293e-05}              & \num{1.376e-04}            & \num{3.075e-06}              \\
      $\frac{1}{80}$       & \num{8.832e-01}            & \num{8.610e-01}              & \num{8.895e-02}            & \num{1.064e-02}              & \num{3.855e-03}            & \num{3.991e-05}              & \num{1.941e-04}            & \num{1.921e-06}              \\
      \bottomrule
    \end{tabular}
  }
\end{table}

\begin{table}[!ht]
  \caption{
  For the $20\times 20$ periodic square inclusion model, the relative errors in the energy norm (in the columns labeled with $\norm{\cdot}_{\tilde{a}}$) and in the $L^2$ norm (in the columns labeled with $\norm{\cdot}_{L^2(\Omega)}$).
  }\label{tab:square-inclusion-error-20}
  \centering
  \makegapedcells
  \footnotesize{
    \begin{tabular}{c c c c c c c c c}
      \toprule
      \multirow{3}{*}{$H$} & \multicolumn{2}{c}{$m=1$}  & \multicolumn{2}{c}{$m=2$}    & \multicolumn{2}{c}{$m=3$}  & \multicolumn{2}{c}{$m=4$}                                                                                                                            \\
      \cmidrule{2-9}
                           & $\norm{\cdot}_{\tilde{a}}$ & $\norm{\cdot}_{L^2(\Omega)}$ & $\norm{\cdot}_{\tilde{a}}$ & $\norm{\cdot}_{L^2(\Omega)}$ & $\norm{\cdot}_{\tilde{a}}$ & $\norm{\cdot}_{L^2(\Omega)}$ & $\norm{\cdot}_{\tilde{a}}$ & $\norm{\cdot}_{L^2(\Omega)}$ \\
      \midrule
      $\frac{1}{10}$       & \num{5.571e-01}            & \num{3.663e-01}              & \num{6.379e-02}            & \num{6.159e-03}              & \num{2.653e-02}            & \num{1.492e-03}              & \num{2.637e-02}            & \num{1.455e-03}              \\
      $\frac{1}{20}$       & \num{4.531e-01}            & \num{2.884e-01}              & \num{3.195e-02}            & \num{1.498e-03}              & \num{2.598e-02}            & \num{1.441e-03}              & \num{2.616e-02}            & \num{1.442e-03}              \\
      $\frac{1}{40}$       & \num{6.306e-01}            & \num{4.938e-01}              & \num{4.544e-02}            & \num{2.613e-03}              & \num{2.589e-02}            & \num{1.361e-03}              & \num{2.782e-02}            & \num{1.485e-03}              \\
      $\frac{1}{80}$       & \num{8.803e-01}            & \num{8.499e-01}              & \num{8.795e-02}            & \num{1.057e-02}              & \num{4.303e-03}            & \num{4.984e-05}              & \num{2.862e-04}            & \num{3.115e-06}              \\
      \bottomrule
    \end{tabular}
  }
\end{table}

\subsection{Periodic cross-shaped inclusion model}
\begin{figure}[!ht]
  \centering
  \includegraphics[width=\textwidth]{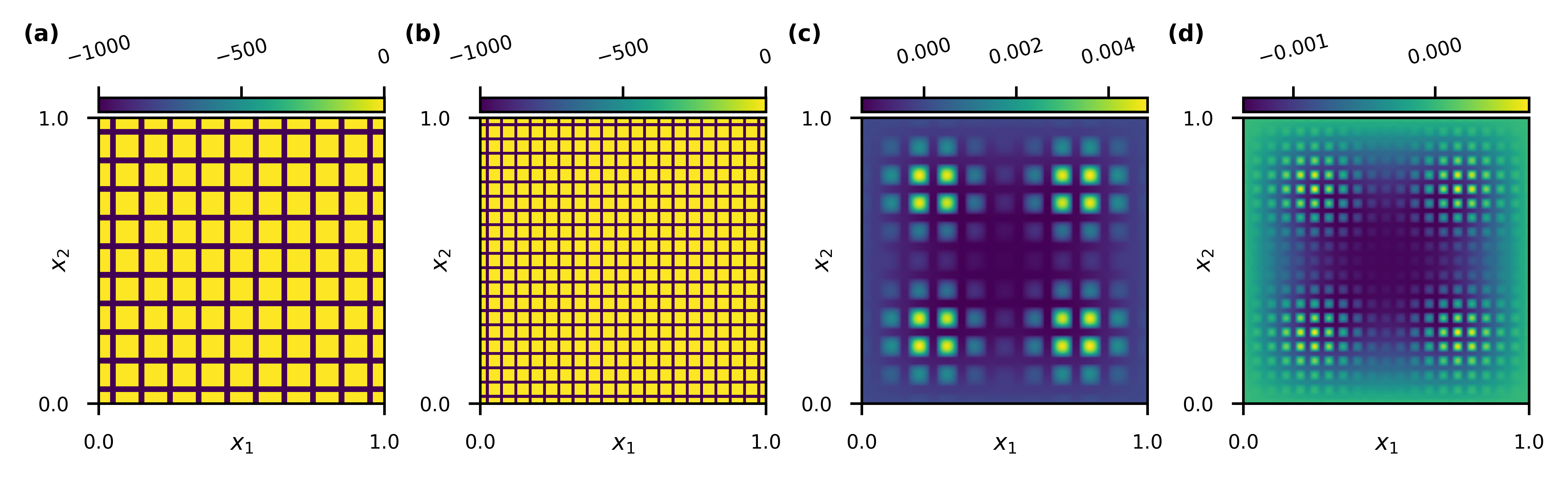}
  \caption{
    Subplots \SubplotTag{(a)} and \SubplotTag{(b)} display visualizations for the coefficient $\sigma$ used in the periodic cross-shaped inclusion model with $10\times 10$ and $20 \times 20$ periodic configurations, respectively.
    Subplots \SubplotTag{(c)} and \SubplotTag{(d)} demonstrate the reference solutions for $10\times 10$ and $20 \times 20$ periodic configurations, respectively.
  }\label{fig:cross-inclusion-sigma1-a}
\end{figure}

In this subsection, we apply the proposed method to a periodic cross-shaped inclusion model.
Specifically, in each periodic cell, a cross-shaped inclusion is centered, imposing a negative coefficient, and the width of the cross arms is set to $1/5$ of the periodic cell.
We consider two periodic configurations, $10\times 10$ and $20\times 20$, as shown in \cref{fig:cross-inclusion-sigma1-a}-\SubplotTag{(a)} and \SubplotTag{(b)}, respectively.
Note that cross-shaped inclusions are all connected, leading to long channels crossing the domain.
Constructing a ``flip'' operator $\mathcal{R}$ in this model to validate the \texttt{T}-coercivity condition is nontrivial compared with the flat interface and the periodic square inclusion models.
Using the theory developed in \cite{BonnetBenDhia2018}, one can prove a weak \texttt{T}-coercivity property holds as long as $\sigma^+_{\mathup{min}}/\sigma^-_{\mathup{max}}$ or $\sigma^-_{\mathup{min}}/\sigma^+_{\mathup{max}}>3$.
This model is designed to test the capability of the proposed method to handle long and high-contrast channels, and we hence set $(\sigma^+_*, \sigma^-_*)=(1, 10^3)$, where $\sigma^+_*$ and $\sigma^-_*$ are defined as in the flat interface model.
The source term $f$ is again from \cref{fig:square-inclusion-source-u}-\SubplotTag{(a)}, and the reference solutions corresponding to the $10\times 10$ and $20\times 20$ periodic configurations are plotted in \cref{fig:cross-inclusion-sigma1-a}-\SubplotTag{(c)} and \SubplotTag{(d)}, respectively.
From reference solution plots, we can observe that long and high-contrast channels create more pronounced oscillations compared to the square inclusion model (cf.\ \cref{fig:square-inclusion-source-u}-\SubplotTag{(b)} and \SubplotTag{(c)}), challenging numerical methods greatly.

\begin{figure}[!ht]
  \centering
  \includegraphics[width=\textwidth]{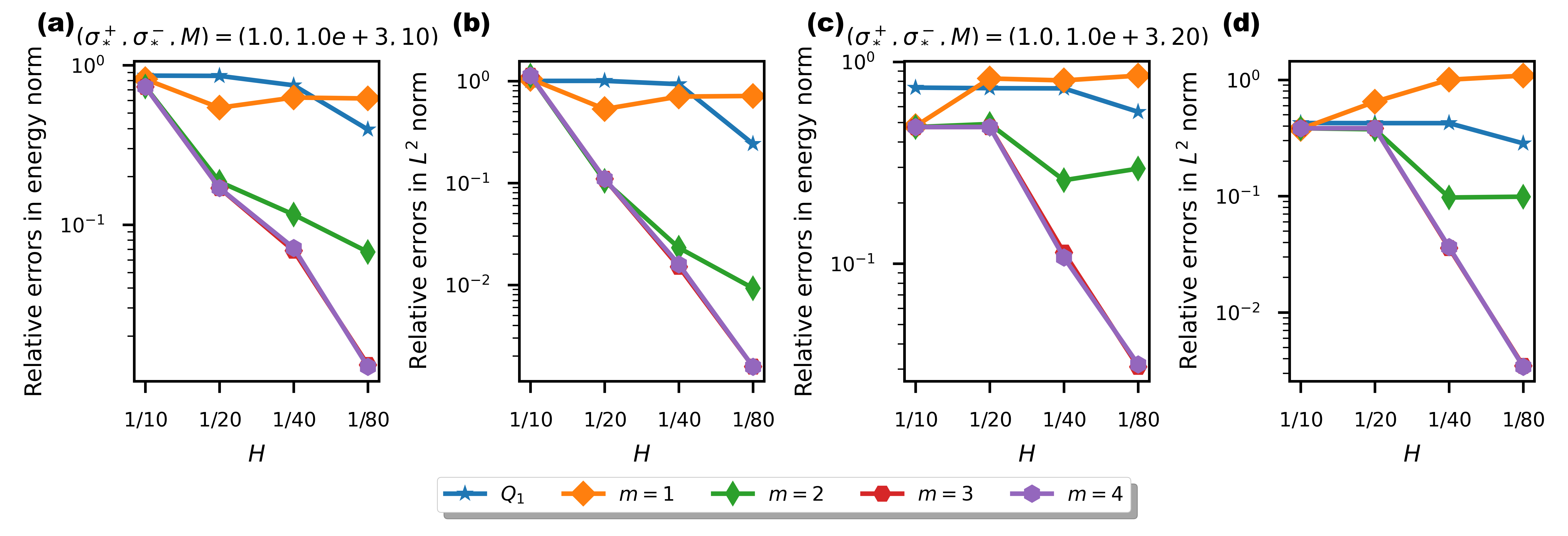}
  \caption{
    Numerical results for the periodic cross-shaped inclusion model with different periodic configurations.
    Subplots \SubplotTag{(a)} and \SubplotTag{(b)} show the relative errors of the proposed method for the $10\times 10$ configuration with different numbers of oversampling layers $m$ and the $Q_1$ FEM w.r.t.\ the coarse mesh size $H$, but measured in different norms.
    Similarly, subplots \SubplotTag{(c)} and \SubplotTag{(d)} correspond to the $20\times 20$ configuration, following the same manner.
  }\label{fig:cross-inclusion-sigma1-b}
\end{figure}

The numerical results of the $Q_1$ FEM and the proposed method with $m \in \{1, 2, 3, 4\}$ are presented in \cref{fig:cross-inclusion-sigma1-b}.
Subplots \SubplotTag{(a)} and \SubplotTag{(b)} in \cref{fig:cross-inclusion-sigma1-b} share the same setting (corresponding to the $10\times 10$ periodic configuration) but the relative errors are measured in the different norms, and similarly for subplots \SubplotTag{(c)} and \SubplotTag{(d)}.
The $Q_1$ FEM fails to provide satisfactory accuracy.
\circled{1} Even with a finer coarse mesh, such as $H=1/80$, there is only a slight improvement, yet the relative errors in the energy norm remain close to $50\%$.
The classical CEM-GMsFEM \cite{Chung2018} is proven to be effective in handling long and high-contrast channels, and the proposed method inherits this advantage.
\circled{2} By setting $m=3$, the proposed method can achieve a relative error of $1\%$ in the energy norm for both the $10\times 10$ and $20\times 20$ periodic configurations.
\circled{3} Typically, the relative errors in the $L^2$ norm are significantly smaller by an order of magnitude than those in the energy norm, and the proposed method can achieve a relative error of $0.1\%$ in the $L^2$ norm for $m=3$ and $4$.
\circled{4} Furthermore, comparing subplots \SubplotTag{(a)} and \SubplotTag{(c)}, as well as \SubplotTag{(b)} and \SubplotTag{(d)}, reveals similar convergence behavior, indicating that the scale of the models does not affect the accuracy of the proposed method.
Note that all coarse meshes employed cannot resolve the channels, highlighting the capability of the proposed method to handle complex coefficient profiles.

\subsection{Random inclusion model}
\begin{figure}[!ht]
  \centering
  \includegraphics[width=0.75\textwidth]{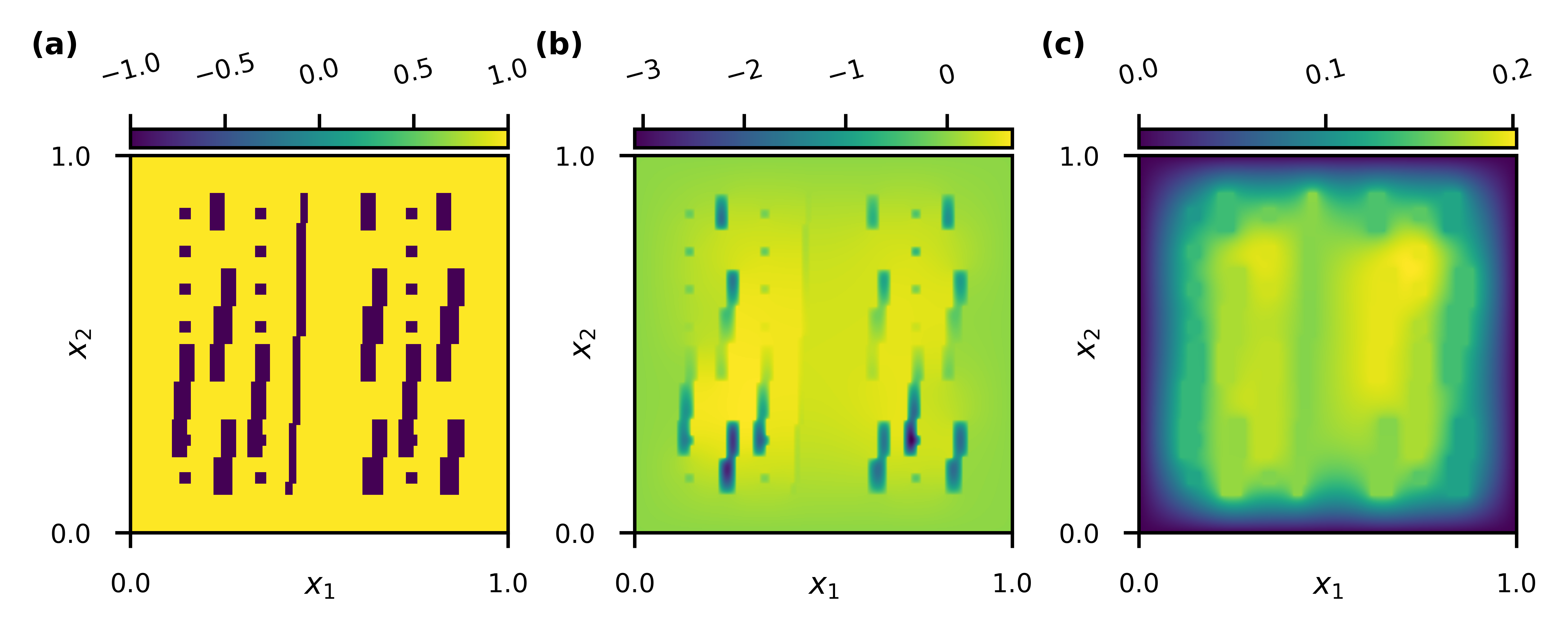}
  \caption{
    \SubplotTag{(a)} The illustration of inclusions (dark regions) of the random inclusion model.
    \SubplotTag{(b)} The plot of the reference solution by setting $(\sigma^+_*, \sigma^-_*)=(1, 10^{-3})$.
    \SubplotTag{(c)} The plot of the reference solution by setting $(\sigma^+_*, \sigma^-_*)=(1, 10^{3})$.
  }\label{fig:random-inclusion-main}
\end{figure}

In this subsection, we consider a random inclusion model that is utilized in several multiscale methods as a showcase of the capability of handling nonperiodic coefficient profiles \cite{Chung2018,Zhao2020,Ye2023a,Poveda2024}.
The subdomains $\Omega^+$ and $\Omega^-$ are demonstrated in \cref{fig:random-inclusion-main}-\SubplotTag{(a)}, and $\sigma$ is determined again by $(\sigma^+_*, \sigma^-_*)$.
We consider two cases: $(\sigma^+_*, \sigma^-_*)=(1, 10^{-3})$ and $(\sigma^+_*, \sigma^-_*)=(1, 10^{3})$.
By setting the source term $f$ as depicted in \cref{fig:square-inclusion-source-u}-\SubplotTag{(a)}, we plot the reference solutions for the two cases, as shown in \cref{fig:random-inclusion-main}-\SubplotTag{(b)} and \SubplotTag{(c)}.
We can observe that void-type inclusions ($(\sigma^+_*, \sigma^-_*)=(1, 10^{-3})$) and rigid-type inclusions ($(\sigma^+_*, \sigma^-_*)=(1, 10^{3})$) exhibit distinct characteristics in the reference solutions.
A recent work \cite{Gorb2021} discussed the asymptotic behavior when the coefficient in inclusions tends to positive infinity, while the case of negative infinity has not been explored in the literature.
Our aim in studying this model is to investigate the effect of choosing different eigenvectors $l^*$.

The approximation of eigenspaces relies on the rapid growth of the eigenvalue $\lambda$ in \cref{eq:new aux}.
Ideally for the Laplace operator, the eigenvalue distribution follows Weyl's law \cite{Weyl1911,Ivrii2016}.
However, as the contrast ratio increases, the growth rate is expected to slow down, leading to a deterioration of the approximation quality.
Surprisingly, this deterioration commonly appears in several leading eigenvalues in \cref{eq:new aux} where a weighted $L^2$ bilinear form is utilized.
To investigate this phenomenon, we examine the first four eigenvalues for $(\sigma^+_*, \sigma^-_*)=(1, 10^{-3})$ and $(\sigma^+_*, \sigma^-_*)=(1, 10^{3})$.
The results are presented in \cref{tab:random-inclusion-sigma-0,tab:random-inclusion-sigma-1}, where the minimal and maximal values of eigenvalues are calculated over all coarse elements.
Since the first eigenvalue $\lambda_1$ is always $0$, we omit the results for $\lambda_1$ in \cref{tab:random-inclusion-sigma-0,tab:random-inclusion-sigma-1}.
\circled{1} We can observe from \cref{tab:random-inclusion-sigma-0} that for $\lambda_2$, there exist small eigenvalues that are on the order of $10^{-4}$.
\circled{2} However, for $\lambda_3$ and $\lambda_4$, the minimum values are significantly larger, on the order of $10^{-1}$, compared to the minimum values of $\lambda_1$.
\circled{3} Interestingly, for $H=1/80$, the maximum and minimum values for both cases are nearly the same.
We conjecture that this is because the coarse meshes are fine enough, and the microstructures on each coarse element are simple.
Note that the eigenvalues by setting $(\sigma^+_*, \sigma^-_*)=(1, 10^{-3})$ and $(\sigma^+_*, \sigma^-_*)=(10^3, 1)$ are strictly identical.
If some coarse elements exhibit a ``symmetrization'' pattern between $\Omega^+$ and $\Omega^-$ and also contribute to the extreme values, we can expect to observe the aforementioned phenomenon.

\begin{table}[!ht]
  \caption{
    For the random inclusion model with $(\sigma^+_*, \sigma^-_*)=(1, 10^{-3})$, the minimal and maximal values of the second/third/fourth eigenvalue over total coarse elements.
  } \label{tab:random-inclusion-sigma-0}
  \centering
  \makegapedcells
  \footnotesize{
    \begin{tabular}{c c c c c c c}
      \toprule
      \multirow{3}{*}{$H$} & \multicolumn{2}{c}{$\lambda_2$} & \multicolumn{2}{c}{$\lambda_3$} & \multicolumn{2}{c}{$\lambda_4$}                                                       \\
      \cmidrule{2-7}
                           & min                             & max                             & min                             & max             & min             & max             \\
      \midrule
      $\frac{1}{10}$       & \num{6.939e-04}                 & \num{4.114e-01}                 & \num{3.011e-01}                 & \num{4.114e-01} & \num{4.121e-01} & \num{8.229e-01} \\
      $\frac{1}{20}$       & \num{7.806e-04}                 & \num{4.121e-01}                 & \num{3.066e-01}                 & \num{8.404e-01} & \num{4.129e-01} & \num{1.278e+00} \\
      $\frac{1}{40}$       & \num{1.008e-03}                 & \num{4.146e-01}                 & \num{2.617e-01}                 & \num{7.579e-01} & \num{5.091e-01} & \num{1.291e+00} \\
      $\frac{1}{80}$       & \num{1.051e-01}                 & \num{6.755e-01}                 & \num{3.255e-01}                 & \num{8.499e-01} & \num{5.228e-01} & \num{1.350e+00} \\
      \bottomrule
    \end{tabular}
  }
\end{table}

\begin{table}[!ht]
  \caption{
    For the random inclusion model with $(\sigma^+_*, \sigma^-_*)=(1, 10^{3})$, the minimal and maximal values of the second/third/fourth eigenvalue over total all coarse elements.
  }\label{tab:random-inclusion-sigma-1}
  \centering
  \makegapedcells
  \footnotesize{
    \begin{tabular}{c c c c c c c}
      \toprule
      \multirow{3}{*}{$H$} & \multicolumn{2}{c}{$\lambda_2$} & \multicolumn{2}{c}{$\lambda_3$} & \multicolumn{2}{c}{$\lambda_4$}                                                       \\
      \cmidrule{2-7}
                           & min                             & max                             & min                             & max             & min             & max             \\
      \midrule
      $\frac{1}{10}$       & \num{2.685e-01}                 & \num{4.335e-01}                 & \num{3.707e-01}                 & \num{1.141e+00} & \num{6.694e-01} & \num{1.553e+00} \\
      $\frac{1}{20}$       & \num{6.856e-02}                 & \num{4.121e-01}                 & \num{3.254e-01}                 & \num{8.696e-01} & \num{5.733e-01} & \num{1.264e+00} \\
      $\frac{1}{40}$       & \num{1.365e-01}                 & \num{6.502e-01}                 & \num{3.292e-01}                 & \num{7.866e-01} & \num{5.763e-01} & \num{1.300e+00} \\
      $\frac{1}{80}$       & \num{1.051e-01}                 & \num{6.755e-01}                 & \num{3.255e-01}                 & \num{8.499e-01} & \num{5.228e-01} & \num{1.350e+00} \\
      \bottomrule
    \end{tabular}
  }
\end{table}

\begin{figure}[!ht]
  \centering
  \includegraphics[width=\textwidth]{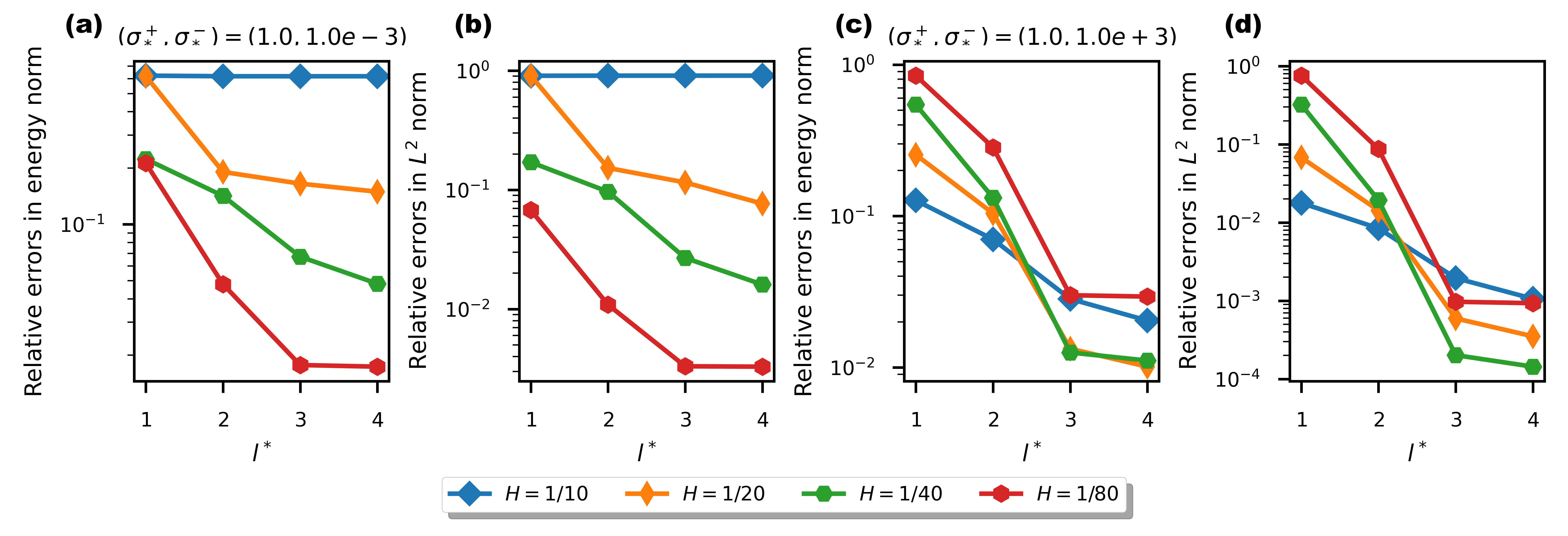}
  \caption{
    Numerical results for the random inclusion model.
    Subplots \SubplotTag{(a)} and \SubplotTag{(b)} show the relative errors of the proposed method for the setting $(\sigma^+_*, \sigma^-_*)=(1, 10^{-3})$ with different coarse mesh size $H$ w.r.t.\ the number of eigenvectors $l^*$, but measured in different norms.
    Similarly, subplots \SubplotTag{(c)} and \SubplotTag{(d)} correspond to $(\sigma^+_*, \sigma^-_*)=(1, 10^{3})$, following the same manner.
  }\label{fig:random-inclusion-error}
\end{figure}

We proceed by presenting the numerical errors of the proposed methods for the two settings in \cref{fig:random-inclusion-error}, where we fix $m=3$ and change $H$ and $l^*$.
Subplots \SubplotTag{(a)} and \SubplotTag{(b)} in \cref{fig:random-inclusion-error} correspond to $(\sigma^+_*, \sigma^-_*)=(1, 10^{-3})$, and subplots \SubplotTag{(c)} and \SubplotTag{(d)} correspond to $(\sigma^+_*, \sigma^-_*)=(1, 10^{3})$.
\circled{1} Upon initial observation, we note a decay of numerical errors with increasing $l^*$ for most cases.
\circled{2} However, this decay is less pronounced with $l^*$ from $3$ to $4$.
Therefore, it is advisable to consider employing a larger number of eigenvectors as a stabilizing strategy rather than solely relying on it for accuracy improvement.
\circled{3} To gain further clarity, the results with $H=1/10$ in \cref{fig:random-inclusion-error}-\SubplotTag{(a)} and \SubplotTag{(b)} do not show a decreasing pattern in errors with $l^*$.
\circled{4} However, we can observe a more substantial reduction in errors as $l^*$ increases for $H=1/10$ in \cref{fig:random-inclusion-error}-\SubplotTag{(c)} and \SubplotTag{(d)}.
This implies that the proposed method is more sensitive to the number of eigenvectors when the coarse mesh is not fine enough, while $l^*=3$ is a recommended choice.

\section{Analysis}\label{sec:analysis}
In this section, we present a rigorous analysis of the proposed method. We begin by examining the global version of the method, where we set $K_i^m$ in \cref{eq:basis} as $\Omega$.
While the global version can yield optimal error estimates, it is impractical due to its high computational cost. Next, we delve into the local version of the method and replicate several estimates found in the original CEM-GMsFEM.
Given the sign-changing setting and the \texttt{T}-coercivity framework, certain assumptions seem to be unavoidable due to technical difficulties, and we will explain the rationale behind these assumptions.

\subsection{Global version}
Let $\Upsilon$ be $\sigma^+_{\mathup{min}}/\sigma^-_{\mathup{max}}$.
As shown in \cref{eq:derivation of T-coercivity}, the well-posedness of the model problem can be ensured if $\Upsilon$ is sufficiently large.
By replacing $\sigma$ with $-\sigma$, the case that $\sigma_{\mathup{min}}^-/\sigma_{\mathup{max}}^+$ is sufficiently large can also be reduced to the case that we mainly focus on.
Essentially, we require that the solution to the model problem be at least unique to discuss the convergence of the numerical method.
Therefore, we introduce the following stronger assumption, although it is rarely touched in the subsequent analysis.
\begin{assumption}\label{ass:zero}
  The model problem \cref{eq:model problem} is well-posed in the sense that there exists a positive constant $C_{\mathup{wp}}$, independent of $f$, such that
  \[
    \norm{u}_{\tilde{a}} \leq C_{\mathup{wp}} \norm{f}_{L^2(\Omega)}.
  \]
\end{assumption}

We define the global operator $\mathcal{G}^\infty_i\colon L^2(\Omega) \rightarrow V=H^1_0(\Omega)$ corresponding to the coarse element $K_i$ via the following variational problem:
\begin{equation}\label{eq:global operator}
  \text{find } \mathcal{G}_i^\infty \psi \in V \text{ s.t. } \forall w \in V,\
  a(\mathcal{G}_i^\infty \psi, w) + s(\mathcal{P}_H\mathcal{G}_i^\infty \psi, \mathcal{P}_H w) = s(\mathcal{P}_H\psi, \mathcal{P}_H w)_{K_i}.
\end{equation}
We reiterate several important facts here: the bilinear form is $s(\cdot, \cdot)$ defined on $L^2(\Omega)\times L^2(\Omega)$ as
\[
  s(v,w) = \int_{\Omega} \mu v w \di x = \mu_{\mathup{msh}} \sum_{i=1}^{N_{\mathup{elem}}} \Diam(K_i)^{-2} \int_{K_i} \sigma v w \di x,
\]
which is not coercive due to the sign-changing property of $\mu$; the operator $\mathcal{P}_H\colon L^2(\Omega)\rightarrow L^2(\Omega)$ is an orthogonal projection with $\im \mathcal{P}_H=V^\mathup{aux}$, regarding the weighted $L^2$ inner-product $\int_\Omega \abs{\mu}v w\di x$ but not the bilinear form $s(\cdot, \cdot)$ or the standard $L^2$ inner-product.
Taking a summation of $\mathcal{G}_i^\infty$ gives $\mathcal{G}^\infty\coloneqq \sum_{i=1}^{N_\mathup{elem}} \mathcal{G}^\infty_i$, which is interpreted as the global operator corresponding to the whole domain.
Certainly, we are required to check the well-posedness of \cref{eq:global operator}, which essentially involves proving the inf-sup stability of the bilinear form:
\[
  a(v, w) + s(\mathcal{P}_Hv, \mathcal{P}_H w),\ \forall (v, w) \in V\times V.
\]
We encounter a challenge here as the properties of $\mathcal{P}_H \mathcal{T}$ are not clear, as $\mathcal{P}_H$ is associated with the coarse mesh $\mathscr{K}_H$, while $\mathcal{T}$ depends on the geometric information of $\Omega^\pm$ and could be rather complicated.
To address this issue, we introduce the following assumption.
\begin{assumption} \label{ass:K_H}
  All coarse elements in $\mathscr{K}_H$ can be categorized into two groups $\mathscr{K}_H^+$ or $\mathscr{K}_H^-$, where
  \[
    \mathscr{K}_H^\pm \coloneqq \CurlyBrackets*{K\mid K \in \mathscr{K}_H \text{ and } K \subset \Omega^\pm}.
  \]
\end{assumption}
This assumption is rather stringent, as it implies that the coarse mesh $\mathscr{K}_H$ is capable of resolving $\Omega^\pm$.
On the other hand, such an assumption greatly facilitates the analysis once the orthogonal projection $\mathcal{T}_H$ comes into play.
Another reason for this assumption is shown in proving \cref{lem:orthogonal}.
We emphasize that several numerical experiments in \cref{sec:numerical experiments} indeed do not align with the assumption, yet the obtained results remain promising.
To facilitate the analysis, we also require the following assumption.
We mention that similar assumptions are also raised in constructing ``flip'' operators that satisfy weak \texttt{T}-coercivity (ref.\ \cite{BonnetBenDhia2018} section 4).
\begin{assumption}\label{ass:L2}
  There exists a positive constant $\norm{\mathcal{R}}_0$, such that for any $v \in V(\Omega^+)$,
  \[
    \norm{\mathcal{R}v}_{0, \Omega^-} \leq \norm{\mathcal{R}}_0 \norm{v}_{0, \Omega^+},
  \]
  where $\mathcal{R}$ is the operator in \cref{eq:T operator}.
\end{assumption}
We introduce this assumption to establish an estimate for $s(\mathcal{P}_Hv, \mathcal{P}_H \mathcal{T}v)$.
This assumption is valid when the operator $\mathcal{R}$ acts as a change of variables.
Denoting $\lambda_i^*$ as the $(l^*+1)$-th eigenvalue of the generalized eigenvalue problem in \cref{eq:new aux}, the following lemma will be frequently utilized in the analysis, while its proof is simply a straightforward application of the properties of eigenspace expansions.
\begin{lemma}\label{lem:eigen}
  On each coarse element $K_i$ and for any $v \in H^1(K_i)$, the following estimates hold:
  \begin{alignat}{1}
    \norm{v-\mathcal{P}_Hv}_{\tilde{s},K_i} & \leq \frac{1}{\sqrt{\lambda_i^*}} \norm{v}_{\tilde{a}, K_i}, \label{eq:eigen a}                                      \\
    \norm{v}_{\tilde{s},K_i}^2              & \leq \norm{\mathcal{P}_Hv}_{\tilde{s},K_i}^2 + \frac{1}{\lambda_i^*} \norm{v}_{\tilde{a}, K_i}^2 \label{eq:eigen b}.
  \end{alignat}
\end{lemma}
We also introduce a notation that $\epsilon \coloneqq 1/\RoundBrackets*{\min_i \lambda_i^*}$.
The following estimates pave the way for proving the well-posedness of \cref{eq:global operator}.
\begin{lemma}
  It holds that for any $v\in V$,
  \begin{align}
    a(v, \mathcal{T}v)                            & \geq \RoundBrackets*{1-\norm{\mathcal{R}}_1/\sqrt{\Upsilon}} \norm{v}_{\tilde{a}}^2, \label{eq:coer a}                                                                                       \\
    s(\mathcal{P}_Hv, \mathcal{P}_H \mathcal{T}v) & \geq \RoundBrackets*{1-\norm{\mathcal{R}}_0/\sqrt{\Upsilon}}\norm{\mathcal{P}_H v}_{\tilde{s}}^2 - \epsilon \norm{\mathcal{R}}_0 / \sqrt{\Upsilon} \norm{v}_{\tilde{a}}^2, \label{eq:coer s} \\
    \norm{\mathcal{T}v}_{\tilde{a}}               & \leq \max\CurlyBrackets*{\RoundBrackets*{1+8\norm{\mathcal{R}}_1^2 / \Upsilon}^{1/2}, \sqrt{2}} \norm{v}_{\tilde{a}}, \label{eq:bound a}                                                     \\
    \norm{\mathcal{T}v}_{\tilde{s}}               & \leq \max\CurlyBrackets*{\RoundBrackets*{1+8\norm{\mathcal{R}}_0^2 / \Upsilon}^{1/2},\sqrt{2}} \norm{v}_{\tilde{s}}. \label{eq:bound s}
  \end{align}
\end{lemma}

\begin{proof}
  The proof of \cref{eq:coer a} has already been given in \cref{eq:derivation of T-coercivity}. We hence first prove \cref{eq:coer s}.
  For any $v\in V$ with $v_1$ and $v_2$ defined as \cref{eq:split of v}, we have
  \begin{align*}
    s(\mathcal{P}_Hv, \mathcal{P}_H \mathcal{T}v) & = \int_{\Omega^+}\abs{\mu}\abs{\mathcal{P}_H v_1}^2 \di x + \int_{\Omega^-}\abs{\mu}\abs{\mathcal{P}_H v_2}^2 \di x-2\int_{\Omega^-}\abs{\mu} \mathcal{P}_Hv_2 \mathcal{P}_H \mathcal{R} v_1 \di x                                                      \\
                                                  & \underset{(\forall \eta > 0)}{\geq} \int_{\Omega^+}\abs{\mu}\abs{\mathcal{P}_H v_1}^2 \di x - \frac{1}{\eta} \int_{\Omega^-} \abs{\mu} \abs{\mathcal{P}_H \mathcal{R} v_1}^2\di x + (1-\eta) \int_{\Omega^-} \abs{\mu} \abs{\mathcal{P}_H v_2}^2 \di x,
  \end{align*}
  where we implicitly utilize \cref{ass:K_H}.
  We turn to the estimate of $\int_{\Omega^-} \abs{\mu} \abs{\mathcal{P}_H \mathcal{R} v_1}^2\di x$ as follows:
  \begin{align*}
    \int_{\Omega^-} \abs{\mu} \abs{\mathcal{P}_H \mathcal{R} v_1}^2\di x & \leq \int_{\Omega^-} \abs{\mu} \abs{\mathcal{R} v_1}^2\di x \quad (\text{\footnotesize $\mathcal{P}_H$ is an orthogonal projection})                                                                                                                                          \\
                                                                         & \leq \mu_{\mathup{max}}^- \int_{\Omega^-} \abs{\mathcal{R}v_1}^2 \di x \leq \mu_{\mathup{max}}^- \norm{\mathcal{R}}_0^2 \int_{\Omega^+} \abs{v_1}^2 \di x \quad (\text{\footnotesize by \cref{ass:L2}})                                                                       \\
                                                                         & \leq \norm{\mathcal{R}}_0^2 \frac{\mu_{\mathup{max}}^-}{\mu_{\mathup{min}}^+} \int_{\Omega^+} \abs{\mu} \abs{v_1}^2 \di x                                                                                                                                                     \\
                                                                         & \leq \norm{\mathcal{R}}_0^2 \frac{\sigma_{\mathup{max}}^-}{\sigma_{\mathup{min}}^+} \RoundBrackets*{\int_{\Omega^+} \abs{\mu} \abs{\mathcal{P}_H v_1}^2\di x+\epsilon \int_{\Omega^+}\abs{\sigma}\abs{\nabla v_1}^2 \di x}. \quad (\text{\footnotesize by \cref{eq:eigen b}})
  \end{align*}
  By choosing $\eta = \norm{\mathcal{R}}_0/\sqrt{\Upsilon}$, we derive that
  \[
    s(\mathcal{P}_Hv, \mathcal{P}_H \mathcal{T}v) \geq (1-\norm{\mathcal{R}}_0/\sqrt{\Upsilon}) \norm{\mathcal{P}_H v}_{\tilde{s}}^2 - \epsilon \norm{\mathcal{R}}_0 / \sqrt{\Upsilon} \norm{v}_{\tilde{a}}^2,
  \]
  which finishes the proof of \cref{eq:coer s}.

  The proofs of \cref{eq:bound a,eq:bound s} follow a similar procedure, and we only provide the proof for the former:
  \begin{align*}
    \norm{\mathcal{T}v}_{\tilde{a}}^2 & = \int_{\Omega^+} \abs{\sigma} \abs{\nabla v_1}^2 \di x + \int_{\Omega^-} \abs{\sigma} \abs{-\nabla v_2+2\nabla\mathcal{R}v_1}^2 \di x                                                                                                        \\
                                      & \leq \int_{\Omega^+} \abs{\sigma} \abs{\nabla v_1}^2 \di x + 2\int_{\Omega^-} \abs{\sigma} \abs{\nabla v_2}^2 \di x + 8 \int_{\Omega^-} \abs{\sigma} \abs{\nabla \mathcal{R}v_1}^2 \di x \quad (\text{\footnotesize by the basic inequality}) \\
                                      & \leq \int_{\Omega^+} \abs{\sigma} \abs{\nabla v_1}^2 \di x + 2\int_{\Omega^-} \abs{\sigma} \abs{\nabla v_2}^2 \di x + 8 \sigma^-_{\mathup{max}} \norm{\mathcal{R}}_{1}^2 \int_{\Omega^+} \abs{\nabla v_1}^2 \di x                             \\
                                      & \leq \RoundBrackets*{1+8\norm{\mathcal{R}}_{1}^2\frac{\sigma^-_{\mathup{max}}}{\sigma^+_{\mathup{min}}}}\int_{\Omega^+} \abs{\sigma} \abs{\nabla v_1}^2 + 2\int_{\Omega^-} \abs{\sigma} \abs{\nabla v_2}^2 \di x                              \\
                                      & \leq \max\CurlyBrackets*{1+8\norm{\mathcal{R}}_1^2 / \Upsilon, 2} \norm{v}_{\tilde{a}}^2.
  \end{align*}
\end{proof}

\begin{remark}
  According to the proof, we can refine $\norm{v}_{\tilde{a}}^2$ in \cref{eq:coer s} to $\norm{v}_{\tilde{a},\Omega^+}^2$.
  However, we choose to retain the original form for the sake of simplicity.
\end{remark}

Now the well-posedness of \cref{eq:global operator} can be established by combining \cref{eq:coer a,eq:coer s}.
\begin{proposition}\label{prop:well-posedness global}
  There exist $\Upsilon'$ and $\epsilon'$ such that for any $\Upsilon \geq \Upsilon'$ and $\epsilon \leq \epsilon'$, the operator $\mathcal{G}_i^\infty$ in \cref{eq:global operator} is well-posed.
\end{proposition}

Recalling that $\im \mathcal{P}_H = V^{\mathup{aux}}$, The following lemma in some sense offers an interpolation operator that maps from $L^2(\Omega)$ to $V$, such that the projections onto $V^\mathup{aux}$ by $\mathcal{P}_H$  are preserved.
\begin{lemma}[ref.\ \cite{Chung2018}]\label{lem:interpolation}
  There exists a bounded map $\mathcal{Q}_H\colon L^2(\Omega) \rightarrow V$ and a positive constant $C_\mathup{inv}$ such that for all $v \in L^2(\Omega)$, it holds that $\mathcal{P}_H\mathcal{Q}_H v=\mathcal{P}_H v$ and $\norm{\mathcal{Q}_H v}_{\tilde{a}} \leq C_\mathup{inv} \norm{\mathcal{P}_H v}_{\tilde{s}}$.
  Moreover, for each coarse element $K_i$, $\mathcal{Q}_Hv|_{K_i}$ depends only on the data of $v$ in $K_i$ and vanishes on $\partial K_i$.
\end{lemma}
We introduce two function spaces as $W\coloneqq \ker \mathcal{P}_H \cap V$ and $V_H^\infty \coloneqq \im \mathcal{G}^\infty \subset V$.
According to \cref{lem:eigen}, it is clear that
\begin{equation}\label{eq:W}
  \norm{w}_{\tilde{s}} \leq \sqrt{\epsilon} \norm{w}_{\tilde{a}},\quad \forall w \in W.
\end{equation}
The following lemma reveals a relationship of ``orthogonality'' between $W$ and $V_H^\infty$ concerning the bilinear form $a(\cdot, \cdot)$.
However, we must be cautious in using the term ``orthogonal'' since $a(\cdot, \cdot)$ cannot define an inner product on $V$.
\begin{lemma}\label{lem:orthogonal}
  For any $v \in V^\infty_H$ and $w \in W$, it holds that $a(v, w)=0$.
  If $w\in V$ and $a(v, w)=0$ holds for any $v \in V^\infty_H$, then $w \in W$.
\end{lemma}
\begin{proof}
  The first argument that $a(v, w)=0$ for any $v \in V^\infty_H$ and $w \in W$ is a direct result of \cref{eq:global operator}.
  We then prove the second argument step by step.
  For simplicity, the notation $\sum_{i,j}$ is used to represent the summation over $i=1,\dots,N_\mathup{elem}$ and $j=1,\dots,l^*$.

  \paragraph{Step1} We first state that the set $\CurlyBrackets{\mathcal{G}^\infty_i\psi_{i,j}\mid 1\leq i \leq N_\mathup{elem}, 1\leq j \leq l^*}$ is linearly independent, where each $\psi_{i,j}$ is an eigenfunction by solving \cref{eq:new aux}.
  Suppose there exists coefficients $\alpha_{i,j}$ such that $\sum_{i,j}\alpha_{i,j}\mathcal{G}^\infty_i\psi_{i,j}=0$.
  Then for any $w\in V$, applying \cref{eq:global operator}, we have
  \[
    0=a(\sum_{i,j}\alpha_{i,j}\mathcal{G}^\infty_i\psi_{i,j}, z) + s(\mathcal{P}_H\sum_{i,j}\alpha_{i,j}\mathcal{G}^\infty_i\psi_{i,j}, \mathcal{P}_H z)=s(\sum_{i,j}\alpha_{i,j}\psi_{i,j}, \mathcal{P}_H z),\quad \forall z \in V.
  \]
  Recalling the property of $\mathcal{Q}_H$ from \cref{lem:interpolation}, we find $s(\sum_{i,j}\alpha_{i,j}\psi_{i,j}, z')=0$ for any $z' \in V^\mathup{aux}$.
  By \cref{ass:K_H}, on each $K_i$, we have $\mu=\abs{\mu}$ or $-\abs{\mu}$, which yields that
  \[
    \int_{K_i} \abs{\mu} \sum_{j} \alpha_{i,j}\psi_{i,j} z'' \di x = \int_{K_i} \mu \sum_{j}  \alpha_{i,j}\psi_{i,j} z'' \di x=0,\quad \forall z'' \in V_i^\mathup{aux}.
  \]
  Therefore, we derive that $\alpha_{i,j}=0$ for any $i$ and $j$, and thus prove the statement.

  \paragraph{Step2} We then state that $\CurlyBrackets{\psi_{i,j}-\mathcal{P}_H \mathcal{G}^\infty_i \psi_{i,j}\mid 1\leq i \leq N_\mathup{elem}, 1\leq j \leq l^*}$ is a linearly independent set.
  Once again, from \cref{eq:global operator}, if $\sum_{i,j}\alpha_{i,j}(\psi_{i,j}-\mathcal{P}_H \mathcal{G}^\infty_i \psi_{i,j})=0$, we have
  \[
    a(\sum_{i,j}\alpha_{i,j}\mathcal{G}^\infty_i\psi_{i,j}, z)=s(\sum_{i,j}\alpha_{i,j}(\psi_{i,j}-\mathcal{P}_H \mathcal{G}^\infty_i \psi_{i,j}), \mathcal{P}_H z)=0,\quad \forall z \in V.
  \]
  According to \cref{ass:zero}, we conclude that $\sum_{i,j}\alpha_{i,j}\mathcal{G}^\infty_i\psi_{i,j}=0$, which yields $\alpha_{i,j}=0$ for any $i$ and $j$ from the previous step.

  \paragraph{Step3} Now we return to the original argument.
  If $\mathcal{P}_H w \neq 0$, we assert that there exists a set of $\CurlyBrackets{\alpha_{i,j}}$ such that
  \[
    s(\sum_{i,j}\alpha_{i,j}(\psi_{i,j}-\mathcal{P}_H \mathcal{G}^\infty_i \psi_{i,j}), \mathcal{P}_H w) = 1.
  \]
  Otherwise, for all choices of $\CurlyBrackets{\alpha_{i,j}}$, we must have $s(\sum_{i,j}\alpha_{i,j}(\psi_{i,j}-\mathcal{P}_H \mathcal{G}^\infty_i \psi_{i,j}), \mathcal{P}_H w) = 0$.
  Taking $\CurlyBrackets{\psi_{i,j}}$ as bases for $V^\mathup{aux}$, we can check that the matrix representation of $s(\cdot, \cdot)$ on $V^\mathup{aux}$ is diagonal, which is also nonsingular thanks to \cref{ass:K_H}.
  From the previous step, we have established that $\CurlyBrackets{\psi_{i,j}-\mathcal{P}_H \mathcal{G}^\infty_i \psi_{i,j}}$ spans the finite-dimentional space $V^\mathup{aux}$.
  However, the relation $s(\sum_{i,j}\alpha_{i,j}(\psi_{i,j}-\mathcal{P}_H \mathcal{G}^\infty_i \psi_{i,j}), \mathcal{P}_H v) = 0$ for any $\CurlyBrackets{\alpha_{i,j}}$ implies that $\mathcal{P}_H v=0$, which is a contradiction.
  Therefore, We have
  \[
    1=a(\sum_{i,j}\alpha_{i,j}\mathcal{G}^\infty_i\psi_{i,j}, w)=s(\sum_{i,j}\alpha_{i,j}(\psi_{i,j}-\mathcal{P}_H \mathcal{G}^\infty_i \psi_{i,j}), \mathcal{P}_H w),
  \]
  which contradicts to $a(v,w)=0$ for any $v \in V^\infty_H$.
\end{proof}

The global solution is defined by solving the following variational form:
\begin{equation}\label{eq:global solution}
  \text{find } u_H^\infty \in V_H^\infty,\ \text{ s.t. }\forall w_H \in V_H^\infty,\ a(u_H, w_H) = \int_\Omega f w_H \di x.
\end{equation}
To guarantee the well-posedness of \cref{eq:global solution}, we need to examine the inf-sup stability of $V_H^\infty$, i.e., proving the existence of a uniform lower bound of
\[
  \inf_{v_H\in V_H^\infty} \sup_{w_H\in V_H^\infty} \frac{a(v_H, w_H)}{\norm{v_H}_{\tilde{a}} \norm{w_H}_{\tilde{a}}}.
\]
The Fortin trick \cite{Boffi2013} suggests that it suffices to check that
\[
  \inf_{v\in W} \sup_{w \in W} \frac{a(v, w)}{\norm{v}_{\tilde{a}} \norm{w}_{\tilde{a}}}\geq \Lambda_W
\]
holds.
The technique that we utilize here is introducing a modified $\mathcal{T}$ operator, denoted as $\mathcal{T}_H$, defined on $W$ as follows:
\begin{align*}
  \mathcal{T}_Hv & = \mathcal{T} v -  \mathcal{Q}_H \mathcal{P}_H \mathcal{T} v=\begin{cases}
v_1-\mathcal{Q}_H\mathcal{P}_H v_1, & \text{ in } \Omega^+, \\
-v_2+2\mathcal{R}v_1 + \mathcal{Q}_H\mathcal{P}_Hv_2 - 2\mathcal{Q}_H \mathcal{P}_H \mathcal{R}v_1, & \text{ in } \Omega^-.
\end{cases} \\
& = \begin{cases}
v_1,    & \text{ in } \Omega^+, \\
-v_2+2\mathcal{R}v_1 - 2\mathcal{Q}_H \mathcal{P}_H \mathcal{R}v_1, & \text{ in } \Omega^-.
\end{cases}
\end{align*}
Here, the terms $\mathcal{Q}_H \mathcal{P}_H v_1$ and $\mathcal{Q}_H \mathcal{P}_H v_2$ vanish because $\mathcal{P}_Hv=0$ on each coarse element, and meanwhile $\mathcal{Q}_H\mathcal{P}_Hv|_{K_i}$ is only dependent on $\mathcal{P}_Hv|_{K_i}$ (see \cref{lem:interpolation}).
Recalling the boundness of $\mathcal{Q}_H$ from \cref{lem:interpolation}, we can see that $\mathcal{T}_H$ is a bounded operator that maps from $W$ to $W$.
Then, we reduce the proof of the inf-sup stability to establish the coercivity of $a(v, \mathcal{T}_Hv)$ for all $v \in W$.
The details of this proof are essentially a replication of \cref{eq:derivation of T-coercivity}.
Taking $\mathcal{R}'_H\coloneqq \mathcal{R}-\mathcal{Q}_H\mathcal{P}_H\mathcal{R}$, we can obtain
\[
  a(v, \mathcal{T}_Hv) \geq \RoundBrackets*{1-\norm{\mathcal{R}'_H}_1/\sqrt{\Upsilon}} \norm{v}_{\tilde{a}}^2,
\]
where $\norm{\mathcal{R}'_H}_1$ is defined similarly as $\norm{\mathcal{R}}_1$ in \cref{eq:norm1 of R}.
Presenting a complete estimate of $\norm{\mathcal{R}'_H}_1$ is complicated.
However, we note that if $\sigma_{\mathup{max}}^-/\sigma_{\mathup{min}}^-$ can be bounded, which is a common scenario in practice, then an apriori estimate of $\norm{\mathcal{R}'_H}_1$ could be achieved.

After proving the inf-sup stability of $V_H^\infty$, an error estimate of the global solution can be derived.
Taking $e=u-u_H^\infty$, we have $a(e, w_H) = 0$ for all $w_H \in V_H^\infty$.
Therefore, \cref{lem:orthogonal} gives that $e\in W$.
Then, we can show that
\begin{align*}
  \norm{e}_{\tilde{a}} & \leq \frac{1}{\Lambda_W} \sup_{w \in W} \frac{a(e, w)}{ \norm{w}_{\tilde{a}}} = \frac{1}{\Lambda_W} \sup_{w \in W} \frac{a(u, w)}{\norm{w}_{\tilde{a}}} = \frac{1}{\Lambda_W} \sup_{w \in W} \frac{\int_\Omega f w \di x}{ \norm{w}_{\tilde{a}}} \\
                       & \leq \frac{1}{\Lambda_W} \sup_{w \in W} \frac{\norm{f}_{\tilde{s}^*} \norm{w}_{\tilde{s}}}{ \norm{w}_{\tilde{a}}}  \quad (\text{\footnotesize by the Cauchy--Schwarz inequality})                                                                \\
                       & \leq \frac{\sqrt{\epsilon}}{\Lambda_W} \norm{f}_{\tilde{s}^*},\quad (\text{\footnotesize by the estimate \cref{eq:W}})
\end{align*}
where
\[
  \norm{f}_{\tilde{s}^*}\coloneqq \RoundBrackets*{\int_\Omega \abs{\mu}^{-1}\abs{f}^2\di x}^{1/2} \approx H \max\CurlyBrackets*{1/\RoundBrackets*{\sigma_\mathup{min}^+}^{1/2}, 1/\RoundBrackets*{\sigma_\mathup{min}^-}^{1/2}} \norm{f}_{L^2(\Omega)}.
\]
We summarize all results presented above in the following theorem.
\begin{theorem}\label{thm:global}
  There exist positive constants $\Upsilon'$ and $\epsilon'$ such that for any $\Upsilon \geq \Upsilon'$ and $\epsilon \leq \epsilon'$, the solution $u_H^\infty$ in \cref{eq:global solution} exists and is unique.
  Moreover, the following estimate
  \[
    \norm{u-u_H^\infty}_{\tilde{a}} \leq \frac{\sqrt{\epsilon}}{\Lambda_W} \norm{f}_{\tilde{s}^*},
  \]
  holds.
\end{theorem}

\subsection{Local version}
When turning to the practical method outlined in \cref{sec:methods}, the initial concern pertains to the existence of multiscale bases as defined in \cref{eq:basis}.
However, a challenge arises due to the incompatibility between the definition of an oversampling region $K_i^m$ and \texttt{T}-coercivity.
Consequently, we cannot obtain a result similar to \cref{prop:well-posedness global} for the local version.
Nevertheless, it is worth noting that the numerical experiments presented in \cref{sec:numerical experiments} did not encounter any lack of well-posedness at the discrete level.
To proceed with the analysis, we introduce the following assumption.
\begin{assumption}\label{ass:modified subdomains}
  There exists a list of subdomains $\CurlyBrackets{\hat{K}_i^1,\hat{K}_i^2,\dots}$ that fulfills the requirements listed below:
  \begin{enumerate}
    \item Each $\hat{K}_i^m$ consists of coarse elements from $\mathscr{K}_H$.
    \item There exists an inclusion relation $K_i \subset \hat{K}_i^1 \subset \hat{K}_i^2 \subset \dots \subset$, such that $\Dist(\partial \hat{K}_i^{m},\partial \hat{K}_i^{m+1})\geq C_\mathup{msh}H$.
    \item For any $v \in V$ with $\Supp v \subset \overline{\hat{K}_i^m}$ or $\overline{\Omega\setminus \hat{K}_i^m}$, it holds that $\Supp \mathcal{T}v \subset \overline{\hat{K}_i^m}$ or $\overline{\Omega\setminus \hat{K}_i^m}$ accordingly.
  \end{enumerate}
\end{assumption}
A construction of such subdomains can be found in \cite{ChaumontFrelet2021}, where a method called symmetrization is described.
However, implementing such a construction (i.e., replacing $K_i^m$ in \cref{eq:basis}) is highly impractical due to the lack of clarity associated with the definition of the operator $\mathcal{T}$.
The fulfillment of the first and second requirements allows us to employ conventional cut-off functions.
Furthermore, the third requirement suggests that, to some extent, \texttt{T}-coercivity can be preserved within these subdomains.
In a sense, this complies with the definition of geometrically-based \texttt{T}-coercivity operators \cite{Nicaise2011,BonnetBenDhia2012,Chesnel2013}.
To distinguish the local problems utilized for analysis from the practical method in \cref{sec:methods}, we add a hat notation ``$\hat{\cdot}$'' (e.g., $\hat{V}_i^m=H_0^1(\hat{K}_i^m)$) to indicate that this term is associated with $\hat{K}_i^m$ rather than $K_i^m$.
Applying the same technique as in establishing \cref{prop:well-posedness global}, we can obtain the well-posedness of the modified local operator $\hat{\mathcal{G}}_i^m$ associated with the coarse element $K_i$:
\begin{equation}\label{eq:local operator}
  \begin{split}
     & \text{find } \hat{\mathcal{G}}_i^m \psi \in \hat{V}_i^m \text{ s.t. } \forall w \in \hat{V}_i^m,                                                                          \\
     & a(\hat{\mathcal{G}}_i^m \psi, w)_{\hat{K}_i^m} + s(\mathcal{P}_H\hat{\mathcal{G}}_i^m \psi, \mathcal{P}_H w)_{\hat{K}_i^m} = s(\mathcal{P}_H\psi, \mathcal{P}_H w)_{K_i}.
  \end{split}
\end{equation}

To facilitate the analysis, we introduce the following notations for constants that are from \cref{eq:coer a,eq:coer s,eq:bound a,eq:bound s}:
\begin{align*}
  C_0 & \coloneqq \max\CurlyBrackets*{\RoundBrackets*{1+8\norm{\mathcal{R}}_0^2 / \Upsilon}^{1/2},\ \sqrt{2}},                                                       \\
  C_1 & \coloneqq \max\CurlyBrackets*{\RoundBrackets*{1+8\norm{\mathcal{R}}_1^2 / \Upsilon}^{1/2},\ \sqrt{2}},                                                       \\
  C_2 & \coloneqq \min\CurlyBrackets*{1-\norm{\mathcal{R}}_1/\sqrt{\Upsilon}-\epsilon\norm{\mathcal{R}}_0/\sqrt{\Upsilon},\ 1-\norm{\mathcal{R}}_0/\sqrt{\Upsilon}}.
\end{align*}
We can see that $C_0$, $C_1$, and $C_2$ can all be bounded from above and below providing that $\Upsilon$ is sufficiently large and $\epsilon$ is sufficiently small, which has also been stated in \cref{prop:well-posedness global,thm:global}.
A key ingredient employed in the analysis is multiplying by a cut-off function $\hat{\chi}_i^{m-1,m}\in C^{0,1}(\Omega)$ defined as
\[
  \hat{\chi}_i^{m-1,m}= \begin{cases}
    1, \text{ in } \hat{K}_i^{m-1}, \\
    0, \text{ in } \Omega\setminus \hat{K}_i^{m},
  \end{cases}
\]
with $0 \leq \hat{\chi}_i^{m-1,m} \leq 1$ in $\hat{K}_i^{m}\setminus \hat{K}_i^{m-1}$. By carefully controlling the parameter $\mu_\mathup{msh}$, we can derive the inequality
\[
  \abs{\nabla \hat{\chi}_i^{m-1,m}}^2 \abs{\sigma} \leq \abs{\mu}
\]
on $\Omega$ for any $i$ and $m$.
Taking $\hat{\mathcal{G}}^m\coloneqq \sum_{i=1}^{N_\mathup{elem}} \hat{\mathcal{G}}^m_i$, We shall prove an estimate of $\mathcal{G}^\infty-\hat{\mathcal{G}}^m$ in \cref{prop:exponential decay} below as the first main result in this subsection, while \cref{lem:prep1,lem:prep2} as presented in \cite{Ye2023a} are reformulated next to reflect the current sign-changing context.

\begin{lemma}\label{lem:prep1}
  There exists a positive constant $\theta$ with $\theta < 1$ that depends on $C_0$, $C_1$, and $C_2$ such that for any $m \geq 1$, $1\leq i \leq N_\mathup{elem}$ and $\psi \in L^2(\Omega)$,
  \[
    \norm{\mathcal{G}_i^\infty\psi}^2_{\tilde{a}, \Omega\setminus \hat{K}_i^m}+\norm{\mathcal{P}_H\mathcal{G}_i^\infty\psi}^2_{\tilde{s}, \Omega\setminus \hat{K}_i^m} \leq \theta^m \RoundBrackets*{\norm{\mathcal{G}_i^\infty\psi}^2_{\tilde{a}}+\norm{\mathcal{P}_H\mathcal{G}_i^\infty\psi}^2_{\tilde{s}}}.
  \]
\end{lemma}

\begin{lemma}\label{lem:prep2}
  It holds that for any $m \geq 1$, $1\leq i \leq N_\mathup{elem}$ and $\psi \in L^2(\Omega)$,
  \[
    \norm{(\mathcal{G}_i^\infty-\hat{\mathcal{G}}_i^m) \psi}_{\tilde{a}}^2 + \norm{\mathcal{P}_H (\mathcal{G}_i^\infty-\hat{\mathcal{G}}_i^m) \psi}_{\tilde{s}}^2 \leq C\theta^{m-1} \RoundBrackets*{\norm{\mathcal{G}_i^\infty\psi}^2_{\tilde{a}}+\norm{\mathcal{P}_H\mathcal{G}_i^\infty\psi}^2_{\tilde{s}}},
  \]
  where $\theta$ here is identical to the one in \cref{lem:prep1}, and $C$ is a positive constant that depends on $C_0$, $C_1$, and $C_2$.
\end{lemma}

\begin{theorem}\label{prop:exponential decay}
  There exist positive constants $\Upsilon'$ and $\epsilon'$ such that for any $\Upsilon \geq \Upsilon'$ and $\epsilon \leq \epsilon'$, it holds that for any $\psi \in L^2(\Omega)$,
  \[
    \norm{(\mathcal{G}^\infty-\hat{\mathcal{G}}^m) \psi}_{\tilde{a}}^2 + \norm{\mathcal{P}_H (\mathcal{G}^\infty-\hat{\mathcal{G}}^m) \psi}_{\tilde{s}}^2 \leq C (m+1)^d \theta^{m-1}\norm{\mathcal{P}_H\psi}_{\tilde{s}}^2,
  \]
  where $\theta$ and $C$ depend on $C_0$, $C_1$, and $C_2$ with $\theta < 1$.
\end{theorem}

In the proof of \cref{prop:exponential decay}, we require an assumption that specifies the growth of the size of $\hat{K}_i^m$ with $m$.
\begin{assumption}\label{ass:overlapping}
  The number of coarse elements within $\hat{K}_i^m$ satisfies a relation:
  \[
    \#\CurlyBrackets*{K\in \mathscr{K}_H\mid K\subset \hat{K}_i^m} \leq C_\mathup{ol}m^d,
  \]
  where $C_\mathup{ol}$ depends solely on the mesh quality.
\end{assumption}

We first prove \cref{lem:prep1}.
\begin{proof}
  Taking $z_i\coloneqq (1-\hat{\chi}_i^{m-1,m})\mathcal{G}_i^\infty \psi$, we can see that $z_i$ is supported in $\Omega\setminus \hat{K}_i^{m-1}$ and hence $(\Supp \mathcal{T}z_i) \cap K_i=\varnothing$ according to \cref{ass:modified subdomains}.3.
  Substituting $\mathcal{T}z_i$ for $w$ in the variational form \cref{eq:global operator}, we have
  \[
    a(\mathcal{G}_i^\infty \psi, \mathcal{T}z_i) + s(\mathcal{P}_H\mathcal{G}_i^\infty \psi, \mathcal{P}_H\mathcal{T}z_i)=0.
  \]
  As a result, we can formulate a decomposition as follows:
  \begin{align*}
     & \quad a(z_i, \mathcal{T}z_i) + s(\mathcal{P}_H z_i, \mathcal{P}_H\mathcal{T}z_i)=-a(\hat{\chi}_i^{m-1,m}\mathcal{G}_i^\infty \psi, \mathcal{T}z_i) - s(\mathcal{P}_H\hat{\chi}_i^{m-1,m}\mathcal{G}_i^\infty \psi, \mathcal{P}_H\mathcal{T}z_i)                                                                                                                                                                             \\
     & = -\underbrace{\int_{\Omega} \sigma \mathcal{G}_i^\infty \psi \nabla \hat{\chi}_i^{m-1,m} \cdot \nabla \mathcal{T}z_i \di x}_{J_1} - \underbrace{\int_{\Omega} \sigma \hat{\chi}_i^{m-1,m} \nabla \mathcal{G}_i^\infty \psi \cdot \nabla \mathcal{T}z_i \di x}_{J_2} - \underbrace{\int_{\Omega} \mu \mathcal{P}_H \RoundBrackets*{\hat{\chi}_i^{m-1,m}\mathcal{G}_i^\infty\psi} \mathcal{P}_H \mathcal{T}z_i \di x}_{J_3}.
  \end{align*}
  For $J_1$, by the Cauchy-Schwarz inequality, it is evident that
  \begin{align*}
    \abs{J_1} & \leq \RoundBrackets*{\int_\Omega \abs{\sigma} \abs{\nabla \hat{\chi}_{i}^{m-1,m}}^2 \abs{\mathcal{G}_i^\infty \psi}^2 \di x }^{1/2} \RoundBrackets*{\int_\Omega \abs{\sigma} \abs{\nabla \mathcal{T} z_i}^2  \di x }^{1/2}                                                                  \\
              & \leq \norm{\mathcal{G}_i^\infty\psi}_{\tilde{s}, \hat{K}_i^{m}\setminus \hat{K}_i^{m-1}} \norm{\mathcal{T}z_i}_{\tilde{a}, \Omega \setminus \hat{K}_i^{m-1}}                                                                                                                                \\
              & \leq C_0\RoundBrackets*{\norm{\mathcal{P}_H\mathcal{G}_i^\infty\psi}_{\tilde{s}, \hat{K}_i^{m}\setminus \hat{K}_i^{m-1}}^2+\epsilon \norm{\mathcal{G}_i^\infty \psi}_{\tilde{a}, \hat{K}_i^{m}\setminus \hat{K}_i^{m-1}}^2}^{1/2} \norm{z_i}_{\tilde{a}, \Omega \setminus \hat{K}_i^{m-1}},
  \end{align*}
  where the last line follows from \cref{eq:bound a} and \cref{lem:eigen}.
  For $J_2$, we have
  \begin{align*}
    \abs{J_2} & \leq \RoundBrackets*{\int_{\hat{K}_i^{m}\setminus \hat{K}_i^{m-1}}\abs{\sigma} \abs{\hat{\chi}_i^{m-1,m}}^2 \abs{\nabla \mathcal{G}_i^\infty \psi}^2 \di x }^{1/2} \RoundBrackets*{\int_\Omega \abs{\sigma} \abs{\nabla \mathcal{T} z_i}^2  \di x }^{1/2} \\
              & \leq \norm{\mathcal{G}_i^\infty\psi}_{\tilde{a}, \hat{K}_i^{m}\setminus \hat{K}_i^{m-1}}\norm{\mathcal{T}z_i}_{\tilde{a}, \Omega \setminus \hat{K}_i^{m-1}}                                                                                               \\
              & \leq C_0\norm{\mathcal{G}_i^\infty\psi}_{\tilde{a}, \hat{K}_i^{m}\setminus \hat{K}_i^{m-1}}\norm{z_i}_{\tilde{a}, \Omega \setminus \hat{K}_i^{m-1}}.
  \end{align*}
  For $J_3$, we can similarly show that
  \begin{align*}
    \abs{J_3} & \leq \RoundBrackets*{\int_{\hat{K}_i^{m}\setminus \hat{K}_i^{m-1}} \abs{\mu} \abs{ \mathcal{P}_H \RoundBrackets*{\hat{\chi}_i^{m-1,m}\mathcal{G}_i^\infty\psi} }^2  \di x }^{1/2} \RoundBrackets*{\int_\Omega \abs{\mu} \abs{\mathcal{P}_H \mathcal{T} z_i}^2  \di x }^{1/2}                                                                                                                          \\
              & = \norm{\mathcal{P}_H \RoundBrackets*{\hat{\chi}_i^{m-1,m}\mathcal{G}_i^\infty\psi}}_{\tilde{s}, \hat{K}_i^{m}\setminus \hat{K}_i^{m-1}} \norm{\mathcal{P}_H \mathcal{T}z_i}_{\tilde{s}, \Omega \setminus \hat{K}_i^{m-1}} \quad (\text{\footnotesize By \cref{ass:modified subdomains}.3})                                                                                                           \\
              & \leq \norm{\hat{\chi}_i^{m-1,m}\mathcal{G}_i^\infty\psi}_{\tilde{s}, \hat{K}_i^{m}\setminus \hat{K}_i^{m-1}} \norm{\mathcal{T}z_i}_{\tilde{s}, \Omega \setminus \hat{K}_i^{m-1}}                                                                                                                                                                                                                      \\
              & \leq C_1\norm{\mathcal{G}_i^\infty\psi}_{\tilde{s}, \hat{K}_i^{m}\setminus \hat{K}_i^{m-1}} \norm{z_i}_{\tilde{s}, \Omega \setminus \hat{K}_i^{m-1}}                                                                                                                                                                                                                                                  \\
              & \leq C_1 \RoundBrackets*{\norm{\mathcal{P}_H\mathcal{G}_i^\infty\psi}_{\tilde{s}, \hat{K}_i^{m}\setminus \hat{K}_i^{m-1}}^2+\epsilon \norm{\mathcal{G}_i^\infty \psi}_{\tilde{a}, \hat{K}_i^{m}\setminus \hat{K}_i^{m-1}}^2}^{1/2} \RoundBrackets*{\norm{\mathcal{P}_Hz_i}_{\tilde{s}, \Omega\setminus \hat{K}_i^{m-1}}^2+\epsilon \norm{z_i}_{\tilde{a}, \Omega \setminus \hat{K}_i^{m-1}}^2}^{1/2}.
  \end{align*}
  Recalling \cref{eq:coer a,eq:coer s}, we can obtain that
  \[
    a(z_i, \mathcal{T}z_i) + s(\mathcal{P}_H z_i, \mathcal{P}_H\mathcal{T}z_i) \geq C_2 \RoundBrackets*{\norm{z_i}_{\tilde{a}, \Omega\setminus \hat{K}_i^{m-1}}^2 + \norm{\mathcal{P}_H z_i}_{\tilde{s}, \Omega\setminus \hat{K}_i^{m-1}}^2}
  \]
  Providing that $\epsilon$ is small enough ($\epsilon \leq 1$), we have
  \begin{align*}
     & \quad \RoundBrackets*{\norm{\mathcal{G}_i^\infty \psi}_{\tilde{a}, \Omega\setminus \hat{K}_i^{m}}^2 + \norm{\mathcal{P}_H \mathcal{G}_i^\infty \psi}_{\tilde{s}, \Omega\setminus \hat{K}_i^{m}}^2}^{1/2} \leq \RoundBrackets*{\norm{z_i}_{\tilde{a}, \Omega\setminus \hat{K}_i^{m-1}}^2 + \norm{\mathcal{P}_H z_i}_{\tilde{s}, \Omega\setminus \hat{K}_i^{m-1}}^2}^{1/2} \\
     & \leq \frac{2C_0+C_1}{C_2} \RoundBrackets*{\norm{\mathcal{G}_i^\infty \psi}_{\tilde{a}, \hat{K}_i^{m}\setminus \hat{K}_i^{m-1}}^2+\norm{\mathcal{P}_H\mathcal{G}_i^\infty\psi}_{\tilde{s}, \hat{K}_i^{m}\setminus \hat{K}_i^{m-1}}^2}^{1/2}.
  \end{align*}
  Note that
  \begin{align*}
     & \quad \norm{\mathcal{G}_i^\infty \psi}_{\tilde{a}, \Omega \setminus \hat{K}_i^{m-1}}^2+\norm{\mathcal{P}_H\mathcal{G}_i^\infty\psi}_{\tilde{s}, \Omega\setminus \hat{K}_i^{m-1}}^2                                                                                                                                                                                    \\
     & =\norm{\mathcal{G}_i^\infty \psi}_{\tilde{a}, \Omega \setminus \hat{K}_i^{m}}^2+\norm{\mathcal{P}_H\mathcal{G}_i^\infty\psi}_{\tilde{s}, \Omega\setminus \hat{K}_i^{m}}^2+\norm{\mathcal{G}_i^\infty \psi}_{\tilde{a}, \hat{K}_i^{m} \setminus \hat{K}_i^{m-1}}^2+\norm{\mathcal{P}_H\mathcal{G}_i^\infty\psi}_{\tilde{s}, \hat{K}_i^{m}\setminus \hat{K}_i^{m-1}}^2.
  \end{align*}
  We hence derive an iteration relation
  \begin{align*}
     & \norm{\mathcal{G}_i^\infty \psi}_{\tilde{a}, \Omega\setminus \hat{K}_i^{m}}^2+\norm{\mathcal{P}_H\mathcal{G}_i^\infty\psi}_{\tilde{s}, \Omega \setminus \hat{K}_i^{m}}^2                                                                                                    \\
     & \qquad \leq  \RoundBrackets*{1+\RoundBrackets*{\frac{C_2}{2C_0+C_1}}^2}^{-1} \RoundBrackets*{\norm{\mathcal{G}_i^\infty \psi}_{\tilde{a}, \Omega \setminus \hat{K}_i^{m-1}}^2+\norm{\mathcal{P}_H\mathcal{G}_i^\infty\psi}_{\tilde{s}, \Omega\setminus \hat{K}_i^{m-1}}^2},
  \end{align*}
  which completes the proof.
\end{proof}

We then prove \cref{lem:prep2}.
\begin{proof}
  We now take $z_i=(\mathcal{G}_i^\infty-\hat{\mathcal{G}}_i^m) \psi$ and introduce a decomposition
  \[
    z_i=\underbrace{(1-\hat{\chi}_i^{m-1,m})\mathcal{G}_i^\infty \psi}_{z'_i} + \underbrace{(\hat{\chi}_i^{m-1,m}-1)\hat{\mathcal{G}}_i^m \psi + \hat{\chi}_i^{m-1,m} z_i}_{z''_i}.
  \]
  We can observe that $z''_i \in \hat{V}_i^m$, which implies $\mathcal{T}z''_i \in \hat{V}_i^m$.
  Therefore, it is easy to see that $a(z_i, \mathcal{T}z''_i)+s(\mathcal{P}_Hz_i, \mathcal{P}_H\mathcal{T}z''_i)=0$ based on the definitions of $\mathcal{G}_i^\infty$ and $\hat{\mathcal{G}}_i^m$ in \cref{eq:global operator,eq:local operator}.
  Now, our task is to estimate $a(z_i, \mathcal{T}z'_i)+s(\mathcal{P}_Hz_i, \mathcal{P}_H\mathcal{T}z'_i)$, and a starting step is
  \[
    a(z_i, \mathcal{T}z'_i)+s(\mathcal{P}_Hz_i, \mathcal{P}_H\mathcal{T}z'_i) \leq \norm{z_i}_{\tilde{a}} \norm{\mathcal{T}z'_i}_{\tilde{a}} + \norm{\mathcal{P}_Hz_i}_{\tilde{s}} \norm{\mathcal{P}_H\mathcal{T}z'_i}_{\tilde{s}}.
  \]
  Recalling that $\norm{\mathcal{T}z'_i}_{\tilde{a}} \leq C_0 \norm{z'_i}_{\tilde{a}}$ and using a similar technique as in the proof of \cref{lem:prep1}, we can obtain that
  \begin{align*}
    \norm{z'_i}_{\tilde{a}}^2 & \leq 2\int_{\Omega \setminus \hat{K}_i^{m-1}}  \abs{\sigma}\abs{1-\hat{\chi}_i^{m-1,m}}^2 \abs{ \nabla {\mathcal{G}}_i^\infty \psi}^2 \di x + 2 \int_{\hat{K}_i^{m}\setminus \hat{K}_i^{m-1}} \abs{\sigma}\abs{\nabla \hat{\chi}_i^{m-1,m}}^2 \abs{{\mathcal{G}}_i^\infty\psi}^2 \di x \\
                              & \leq 2 \norm{{\mathcal{G}}_i^\infty \psi}^2_{\tilde{a},\Omega \setminus \hat{K}_i^{m-1}} + 2 \int_{\Omega\setminus \hat{K}_i^{m-1}} \abs{\mu} \abs{\mathcal{G}_i^\infty \psi}^2 \di x \quad (\text{\footnotesize by \cref{ass:modified subdomains}.2})                                 \\
                              & \leq 2(1+\epsilon)\norm{{\mathcal{G}}_i^\infty \psi}^2_{\tilde{a},\Omega \setminus \hat{K}_i^{m-1}} + 4 \norm{\mathcal{P}_H \mathcal{G}_i^\infty \psi}_{\tilde{s}, \Omega \setminus \hat{K}_i^{m-1}}^2.
  \end{align*}
  Moreover, we can see that $\norm{\mathcal{P}_H\mathcal{T}z'_i}_{\tilde{s}} \leq \norm{\mathcal{T}z'_i}_{\tilde{s}} \leq C_1\norm{z'_i}_{\tilde{s}} $ and
  \begin{align*}
    \norm{z'_i}_{\tilde{s}}^2 & =\int_{\Omega} \abs{\mu} \abs{1-\hat{\chi}_i^{m-1,m}}^2 \abs{\mathcal{G}_i^\infty \psi}^2 \di x \leq \norm{\mathcal{G}_i^\infty \psi}_{\tilde{s},\Omega\setminus \hat{K}_i^{m-1}}^2          \\
                              & \leq \epsilon \norm{\mathcal{G}_i^\infty \psi}_{\tilde{a},\Omega\setminus \hat{K}_i^{m-1}}^2 + \norm{\mathcal{P}_H\mathcal{G}_i^\infty \psi}_{\tilde{s}, \Omega\setminus \hat{K}_i^{m-1}}^2.
  \end{align*}
  Then, we are close to the target estimate as
  \begin{align*}
    C_2\RoundBrackets*{\norm{z_i}_{\tilde{a}}^2 + \norm{\mathcal{P}_H z_i}_{\tilde{s}}^2} & \leq a(z_i, \mathcal{T}z_i) + s(\mathcal{P}_Hz_i, \mathcal{P}_H \mathcal{T}z_i)=a(z_i, \mathcal{T}z'_i)+s(\mathcal{P}_Hz_i, \mathcal{P}_H\mathcal{T}z'_i)                                                    \\
                                                                                          & \leq \RoundBrackets*{\norm{z_i}_{\tilde{a}}^2 + \norm{\mathcal{P}_H z_i}_{\tilde{s}}^2}^{1/2} \RoundBrackets*{\norm{\mathcal{T}z'_i}_{\tilde{a}}^2+\norm{\mathcal{P}_H\mathcal{T}z'_i}_{\tilde{s}}^2}^{1/2}.
  \end{align*}
  Suppose that $\epsilon$ is small enough, we can obtain that
  \[
    \norm{z_i}_{\tilde{a}}^2 + \norm{\mathcal{P}_H z_i}_{\tilde{s}}^2 \leq \RoundBrackets*{\frac{4C_0+C_1}{C_2}}^2 \RoundBrackets*{\norm{\mathcal{G}_i^\infty \psi}^2_{\tilde{a},\Omega\setminus \hat{K}_i^{m-1}}+\norm{\mathcal{P}_H \mathcal{G}_i^\infty \psi}_{\tilde{s},\Omega\setminus \hat{K}_i^{m-1}}^2},
  \]
  and we hence derive the desired result by utilizing \cref{lem:prep1}.
\end{proof}

We are now ready to prove \cref{prop:exponential decay}.
\begin{proof}
  We denote $z_j\coloneqq (\mathcal{G}_j^\infty-\hat{\mathcal{G}}_j^m)\psi$ and $z \coloneqq \sum_{j=1}^{N_\mathup{elem}} z_j$.
  We take a decomposition of $z$ as
  \[
    z = \underbrace{(1-\hat{\chi}_i^{m,m+1})z_i}_{z'}+\underbrace{\hat{\chi}_i^{m,m+1}z}_{z''},
  \]
  where $i$ is arbitrary chosen from $1,\dots,N_\mathup{elem}$.
  Again, we can observe that $\Supp(\mathcal{T}z') \cap K_i=\varnothing$ and hence have
  \[
    a(\mathcal{G}_i^\infty\psi, \mathcal{T}z')+s(\mathcal{P}_H\mathcal{G}_i^\infty\psi,\mathcal{P}_H\mathcal{T}z') = a(\hat{\mathcal{G}}_i^m\psi, \mathcal{T}z')+s(\mathcal{P}_H\hat{\mathcal{G}}_i^m\psi,\mathcal{P}_H\mathcal{T}z')=0,
  \]
  which leads to
  \[
    a(z_i, \mathcal{T}z)+s(\mathcal{P}_Hz_i,\mathcal{P}_H\mathcal{T}z)=a(z_i, \mathcal{T}z'')+s(\mathcal{P}_Hz_i,\mathcal{P}_H\mathcal{T}z'')\leq \norm{\mathcal{T}z''}_{\tilde{a}}\norm{z_i}_{\tilde{a}}+\norm{\mathcal{P}_H\mathcal{T}z''}_{\tilde{s}}\norm{\mathcal{P}_Hz_i}_{\tilde{s}}.
  \]
  Similarly, we can estimate $\norm{\mathcal{T}z''}_{\tilde{a}}$ and $\norm{\mathcal{P}_H\mathcal{T}z''}_{\tilde{s}}$ as
  \[
    \norm{\mathcal{T}z''}_{\tilde{a}} \leq C_0 \norm{\hat{\chi}_i^{m,m+1}z}_{\tilde{a}}\leq C_0\RoundBrackets*{\norm{z}_{\tilde{a},\hat{K}_i^{m+1}}+\norm{z}_{\tilde{s},\hat{K}_i^{m+1}}}\leq C_0\RoundBrackets*{(1+\sqrt{\epsilon}) \norm{z}_{\tilde{a},\hat{K}_i^{m+1}} + \norm{\mathcal{P}_Hz}_{\tilde{s},\hat{K}_i^{m+1}}}
  \]
  and
  \[
    \norm{\mathcal{P}_H\mathcal{T}z''}_{\tilde{s}}\leq \norm{\mathcal{T}z''}_{\tilde{s}} \leq C_1 \norm{\hat{\chi}_i^{m,m+1}z}_{\tilde{s}} \leq C_1\norm{z}_{\tilde{s},\hat{K}_i^{m+1}} \leq C_1\RoundBrackets*{\sqrt{\epsilon}\norm{z}_{\tilde{a},\hat{K}_i^{m+1}}+\norm{\mathcal{P}_Hz}_{\tilde{s},\hat{K}_i^{m+1}}}.
  \]
  To simplify the expression, we assume that $\epsilon \leq 1$ and derive that
  \[
    a(z_i, \mathcal{T}z)+s(\mathcal{P}_Hz_i,\mathcal{P}_H\mathcal{T}z)\leq c_0\RoundBrackets*{\norm{z}_{\tilde{a},\hat{K}_i^{m+1}}^2+\norm{\mathcal{P}_Hz}_{\tilde{s},\hat{K}_i^{m+1}}^2}^{1/2} \RoundBrackets*{\norm{z_i}_{\tilde{a}}^2+\norm{\mathcal{P}_H z_i}_{\tilde{s}}^2}^{1/2},
  \]
  where $c_0$ depends on $C_0$, $C_1$, and $C_2$.
  Based on \cref{ass:overlapping}, it holds that
  \[
    \sum_{i=1}^{N_\mathup{elem}} \RoundBrackets*{\norm{z}_{\tilde{a},\hat{K}_i^{m+1}}^2+\norm{\mathcal{P}_Hz}_{\tilde{s},\hat{K}_i^{m+1}}^2} \leq C_\mathup{ol}(m+1)^d \RoundBrackets*{ \norm{z}_{\tilde{a}}^2+\norm{\mathcal{P}_Hz}_{\tilde{s}}^2}.
  \]
  Therefore, applying the Cauchy--Schwarz inequality,
  \begin{align*}
     & \quad C_2 \RoundBrackets*{\norm{z}_{\tilde{a}}^2+\norm{\mathcal{P}_Hz}_{\tilde{s}}^2} \leq \sum_{i=1}^{N_\mathup{elem}} \RoundBrackets*{a(z_i, \mathcal{T}z)+s(\mathcal{P}_Hz_i,\mathcal{P}_H\mathcal{T}z)}                                                                                                  \\
     & \leq c_0 \RoundBrackets*{\sum_{i=1}^{N_\mathup{elem}} \RoundBrackets*{\norm{z}_{\tilde{a},\hat{K}_i^{m+1}}^2+\norm{\mathcal{P}_Hz}_{\tilde{s},\hat{K}_i^{m+1}}^2}}^{1/2} \RoundBrackets*{\sum_{i=1}^{N_\mathup{elem}}\RoundBrackets*{\norm{z_i}_{\tilde{a}}^2+\norm{\mathcal{P}_H z_i}_{\tilde{s}}^2}}^{1/2} \\
     & \leq c_0 \sqrt{C_\mathup{ol}} (m+1)^{d/2}\RoundBrackets*{ \norm{z}_{\tilde{a}}^2+\norm{\mathcal{P}_Hz}_{\tilde{s}}^2}^{1/2} \RoundBrackets*{\sum_{i=1}^{N_\mathup{elem}}\RoundBrackets*{\norm{z_i}_{\tilde{a}}^2+\norm{\mathcal{P}_H z_i}_{\tilde{s}}^2}}^{1/2}.
  \end{align*}
  It has been shown in \cref{lem:prep2} that, for any $i$,
  \[
    \norm{z_i}_{\tilde{a}}^2+\norm{\mathcal{P}_H z_i}_{\tilde{s}}^2 \leq  C_*\theta^{m-1}\RoundBrackets*{\norm{\mathcal{G}_i^\infty \psi}_{\tilde{a}}^2 + \norm{\mathcal{P}_H\mathcal{G}_i^\infty \psi}_{\tilde{s}}^2},
  \]
  and we turn to provide a bound on $\sum_{i=1}^{N_\mathup{elem}} \norm{\mathcal{G}_i^\infty \psi}_{\tilde{a}}^2 + \norm{\mathcal{P}_H\mathcal{G}_i^\infty \psi}_{\tilde{s}}^2$.
  Note that by the definition of \cref{eq:local operator},
  \begin{align*}
     & \quad  C_2 \CurlyBrackets*{\norm{\mathcal{G}_i^\infty \psi}_{\tilde{a}}^2 + \norm{\mathcal{P}_H\mathcal{G}_i^\infty \psi}_{\tilde{s}}^2} \leq a(\mathcal{G}_i^\infty\psi, \mathcal{T}\mathcal{G}_i^\infty\psi)+s(\mathcal{P}_H\mathcal{G}_i^\infty\psi, \mathcal{P}_H \mathcal{T} \mathcal{G}_i^\infty\psi)   \\
     & = s(\mathcal{P}_H\psi, \mathcal{P}_H \mathcal{T}\mathcal{G}_i^\infty \psi)_{K_i} \leq \norm{\mathcal{P}_H\psi}_{\tilde{s},K_i}\norm{\mathcal{P}_H\mathcal{T}\mathcal{G}_i^\infty\psi}_{\tilde{s},K_i} \leq \norm{\mathcal{P}_H\psi}_{\tilde{s},K_i}\norm{\mathcal{T}\mathcal{G}_i^\infty\psi}_{\tilde{s},K_i} \\
     & \leq \norm{\mathcal{P}_H\psi}_{\tilde{s},K_i}\norm{\mathcal{T}\mathcal{G}_i^\infty\psi}_{\tilde{s}} \leq C_1\norm{\mathcal{P}_H\psi}_{\tilde{s},K_i}\norm{\mathcal{G}_i^\infty\psi}_{\tilde{s}}                                                                                                               \\
     & \leq C_1\norm{\mathcal{P}_H\psi}_{\tilde{s},K_i}\RoundBrackets*{\sqrt{\epsilon}\norm{\mathcal{G}_i^\infty\psi}_{\tilde{a}}+\norm{\mathcal{P}_H\mathcal{G}_i^\infty\psi}_{\tilde{s}}},
  \end{align*}
  which gives that
  \[
    \norm{\mathcal{G}_i^\infty \psi}_{\tilde{a}}^2 + \norm{\mathcal{P}_H\mathcal{G}_i^\infty \psi}_{\tilde{s}}^2 \leq c_1 \norm{\mathcal{P}_H\psi}_{\tilde{s},K_i}^2.
  \]
  where $c_1$ depends on $C_1$ and $C_2$. We finally completed the proof by collecting all the estimates obtained above.
\end{proof}

We similarly denote $\hat{\mathcal{G}}^m \coloneqq \sum_{i=1}^{N_\mathup{elem}} \hat{\mathcal{G}}^m_i$, where we implicitly lift the image of $\hat{\mathcal{G}}^m_i$ to $V$.
We now define the multiscale solution in the local version as the solution to the following problem: for $m \geq 1$,
\begin{equation}\label{eq:local ms}
  \text{find } \hat{u}_H^m \in \hat{V}^m_H \text{ s.t. }\forall w_H \in \hat{V}^m_H,\ a(\hat{u}_H^m, w_H) = \int_{\Omega} f w_H \di x,
\end{equation}
where $\hat{V}^m_H \coloneqq \im \hat{\mathcal{G}}^m \subset V$.
The remainder of this section is dedicated to proving an error estimate for the multiscale solution, namely \cref{thm:main local} below.
Once again, we need to verify the well-posedness of \cref{eq:local ms}, or equivalently, demonstrate that
\[
  \inf_{v_H\in \hat{V}_H^m} \sup_{w_H\in \hat{V}_H^m} \frac{a(v_H, w_H)}{\norm{v_H}_{\tilde{a}} \norm{w_H}_{\tilde{a}}}
\]
can be bounded from below by a positive constant that is independent of $m$ and $H$.
Thanks to \cref{thm:global}, the inf-sup condition holds for the global problem, while a global basis and a local basis are connected through $\psi$.
Those are the main ingredients to prove the inf-sup condition for the local problem. First, for any $v_H \in \hat{V}_H^m$, we can find $\psi$ such that $v_H = \hat{\mathcal{G}}^m \psi$.
Next, We choose $v_H' = \mathcal{G}^\infty \psi \in V_H^\infty$.
Recalling the inf-sup stability on $V^\infty_H$, we can find $w_H' = \mathcal{G}^\infty \phi \in V_H^\infty$ such that $a(v_H', w_H') \geq \Lambda^\infty \norm{v_H'}_{\tilde{a}} \norm{w_H'}_{\tilde{a}}$.
Similarly, if we denote $w_H = \hat{\mathcal{G}}^m \phi$, we have
\begin{align*}
  a(v_H, w_H) & = a(v_H',w_H')+a(v_H',w_H-w_H')+a(v_H-v_H',w_H)                                                                                                                           \\
              & \geq \Lambda^\infty \norm{v_H'}_{\tilde{a}}\norm{w_H'}_{\tilde{a}}-\norm{v_H'}_{\tilde{a}}\norm{w_H-w_H'}_{\tilde{a}}-\norm{v_H-v_H'}_{\tilde{a}} \norm{w_H}_{\tilde{a}}.
\end{align*}
Therefore, if we can show that
\begin{equation}\label{eq:target}
  \norm{(\mathcal{G}^\infty-\hat{\mathcal{G}}^m)\psi}_{\tilde{a}} \leq \eta \norm{\mathcal{G}^\infty \psi}_{\tilde{a}},
\end{equation}
where $\eta$ is a small constant, we can finish the proof.
The reason is that we now have $ 1/2 \norm{v_H}_{\tilde{a}} \leq \norm{v'_H}_{\tilde{a}} \leq 3/2 \norm{v_H}_{\tilde{a}}$ and $1/2 \norm{w_H}_{\tilde{a}} \leq \norm{w'_H}_{\tilde{a}} \leq 3/2 \norm{w_H}_{\tilde{a}}$ by the smallness of $\eta$.
Consequently, we can deduce that
\[
  a(v_H, w_H) \geq \RoundBrackets*{{\Lambda^\infty}/{4}-{9\eta}/{4} - {3\eta}/{2}} \norm{v_H}_{\tilde{a}} \norm{w_H}_{\tilde{a}}.
\]
By comparing \cref{eq:target} with \cref{prop:exponential decay}, our remaining task is to demonstrate $\norm{\mathcal{P}_H \psi}_{\tilde{s}}$ can be bounded by $\norm{\mathcal{G}^\infty \psi}_{\tilde{a}}$.
According to \cref{lem:interpolation}, it not harmful to assume that $\psi \in V$ with $\norm{\psi}_{\tilde{a}} \leq C_\mathup{inv} \norm{\mathcal{P}_H\psi}_{\tilde{s}}$.
Taking $\mathcal{T}\psi$ as a test function, we can obtain that
\[
  a(\mathcal{G}^\infty\psi, \mathcal{T} \psi)+ s(\mathcal{P}_H\mathcal{G}^\infty\psi, \mathcal{P}_H\mathcal{T}\psi)=s(\mathcal{P}_H\psi, \mathcal{P}_H\mathcal{T} \psi).
\]
The estimate \cref{eq:coer s} says that
\[
  C'\norm{\mathcal{P}_H \psi}_{\tilde{s}}^2-C''\norm{\psi}_{\tilde{a}}^2 \leq s(\mathcal{P}_H\psi, \mathcal{P}_H\mathcal{T} \psi)
\]
where $C' \rightarrow 1$ and $C'' \rightarrow 0$ as $\Upsilon \rightarrow \infty$.
Therefore, if $\Upsilon$ is large enough, we can show that
\begin{align*}
  c_0\norm{\mathcal{P}_H \psi}_{\tilde{s}}^2 & \leq  \RoundBrackets*{C'-C''C_\mathup{inv}^2}\norm{\mathcal{P}_H \psi}_{\tilde{s}}^2 \leq C'\norm{\mathcal{P}_H \psi}_{\tilde{s}}^2-C''\norm{\psi}_{\tilde{a}}^2                                                                                     \\
                                             & \leq s(\mathcal{P}_H\psi, \mathcal{P}_H\mathcal{T} \psi) \leq \norm{\mathcal{G}^\infty\psi}_{\tilde{a}} \norm{\mathcal{T}\psi}_{\tilde{a}}+\norm{\mathcal{P}_H\mathcal{G}^\infty\psi}_{\tilde{s}} \norm{\mathcal{P}_H\mathcal{T}\psi}_{\tilde{s}}    \\
                                             & \leq c_1\sqrt{1+C_\mathup{inv}^2}\RoundBrackets*{\norm{\mathcal{G}^\infty\psi}_{\tilde{a}}^2+\norm{\mathcal{P}_H\mathcal{G}^\infty\psi}_{\tilde{s}}^2}^{1/2}\RoundBrackets*{\norm{\psi}_{\tilde{a}}^2+\norm{\mathcal{P}_H\psi}_{\tilde{s}}^2}^{1/2},
\end{align*}
where $c_0$ and $c_1$ are positive constants and the last line follows the techniques in the previous proofs.
We hence obtain that
\[
  \norm{\psi}_{\tilde{a}}/C_\mathup{inv} \leq \norm{\mathcal{P}_H\psi}_{\tilde{s}} \leq c_2 \RoundBrackets*{\norm{\mathcal{G}^\infty\psi}_{\tilde{a}}^2+\norm{\mathcal{P}_H\mathcal{G}^\infty\psi}_{\tilde{s}}^2}^{1/2}.
\]
Meanwhile, utilizing Poincar\'e inequality, we can see that
\[
  \norm{\mathcal{P}_H\mathcal{G}^\infty\psi}_{\tilde{s}} \leq \norm{\mathcal{G}^\infty\psi}_{\tilde{s}} \leq C_\mathup{po} H^{-1} \norm{\mathcal{G}^\infty\psi}_{\tilde{a}}.
\]
Hence, in conjunction with \cref{prop:exponential decay}, we establish the existence of such $\eta$ in \cref{eq:target}.
Once the inf-sup condition is validated, the existence and uniqueness of the multiscale solution are guaranteed.
Moreover, we can employ C\'{e}a's lemma to derive the error estimate.
Specifically, we obtain
\[
  \norm{u-\hat{u}_H^m}_{\tilde{a}} \leq c_3\norm{u-v_H}_{\tilde{a}} \leq c_3\RoundBrackets*{\norm{u-u_H^\infty}_{\tilde{a}}+\norm{u_H^\infty-v_H}_{\tilde{a}}},
\]
where $v_H=\hat{\mathcal{G}}^m\psi \in \hat{V}^m_H$ if $u_H^\infty=\mathcal{G}^\infty \psi$.
We have already shown that
\[
  \norm{u_H^\infty-v_H}_{\tilde{a}} \leq c_4 H^{-1}(m+1)^{d/2}\theta^{(m-1)/2} \norm{u_H^\infty}_{\tilde{a}}.
\]
Finally, we present the main theorem of this section.

\begin{theorem}\label{thm:main local}
  There exist positive constants $\Upsilon'$ and $\epsilon'$ such that for any $\Upsilon \geq \Upsilon'$ and $\epsilon \leq \epsilon'$, the solution to \cref{eq:local ms} exists and is unique.
  Moreover, the following error estimate holds:
  \[
    \norm{u-\hat{u}_H^m}_{\tilde{a}} \leq C_* \RoundBrackets*{1+ H^{-2}(m+1)^{d/2}\theta^{(m-1)/2}} \norm{f}_{\tilde{s}^*},
  \]
  where the positive constant $C_*$ and $\theta$ is independent of $H$ and $m$ with $\theta < 1$.
\end{theorem}

\section{Conclusions}\label{sec:conclusions}
We have proposed a multiscale computational method for solving sign-changing problems, utilizing the framework of CEM-GMsFEMs in the construction of multiscale basis functions.
However, a direct application of the original CEM-GMsFEM encounters an immediate challenge during the construction of auxiliary spaces, as the generalized spectral problems can become ill-defined due to a non-positive definite right-hand bilinear form.
We have addressed this issue by replacing the coefficient with its absolute value in this step and explained that this modification complies with the \texttt{T}-coercivity theory.
Moreover, we have focused on the relaxed version of the CEM-GMsFEM, which is more implementation-friendly as it eliminates the need to solve saddle-point problems.
The numerical experiments conducted have highlighted several advantages of the proposed method: (1) the flexibility that coarse meshes do not require to resolve with interfaces, (2) the accuracy that remains stable under low regular exact solutions, and (3) the robustness in high contrast coefficient profiles.
Due to technical difficulties, such as the nonlocality of the reflection operator $\mathcal{T}$, the final error estimates have been proved under several stringent assumptions, primarily related to coarse element partitions and coefficient contrast.
However, we emphasize that the potential of the proposed method has been demonstrated through numerical experiments, and we firmly believe that the theoretical results can be further improved by relaxing the assumptions and developing more sophisticated analysis techniques.

The current implementation of the proposed method did not exploit the parallelism in constructing multiscale basis functions, and hence the computation time in the offline stage is awkwardly long, even longer than generating reference solutions on the fine mesh.
Therefore, we plan to investigate the parallelization in the shared memory architecture of the offline stage in the future.
We emphasize that the application of computational multiscale methods to sign-changing problems is much more appealing in large-scale simulations.
As a matter of fact, the algebraic linear systems in this situation are no longer positive definite, which renders iterative solvers less efficient, even with preconditioning techniques.
Consequently, in this context of direct solvers, the reduction of degrees of freedom by multiscale methods is expected to be more beneficial.

\noindent {\bf Acknowledgments.} EC's research is partially supported by the Hong Kong RGC General Research Fund (Project numbers: 14305624 and 14305423).


\end{document}